\documentclass[twoside,11pt,leqno]{amsart}
\usepackage[latin1]{inputenc}
\usepackage{anysize}
\usepackage{amsfonts,amsopn,amsthm,amsmath,amssymb,latexsym}
\newcommand{\cal}{\mathcal}
\newcommand{\C}{\mathbb{C}}
\newcommand{\R}{\mathbb{R}}
\newcommand{\Z}{\mathbb{Z}}
\newcommand{\K}{\mathbb{K}}
\newcommand{\ff}{{\bf f}}
\newcommand{\fg}{{\bf g}}
\newcommand{\rH}{{\rm H}}
\newcommand{\cA}{{\cal A}}

\newcommand{\cBS}{{\cal F^{\BS}}}

\newcommand{\cF}{{\cal F}}

\newcommand{\id}{{\text{\rm id}}}

\newcommand{\gmin}{{\gamma_{\rm min}}}
\newcommand{\gmax}{{\gamma_{\rm max}}}
\DeclareMathOperator{\im}{im}

\DeclareMathOperator{\Ht}{ht}
\DeclareMathOperator{\Eu}{Eu}

\DeclareMathOperator{\Span}{span}
\DeclareMathOperator{\Sym}{Sym}
\DeclareMathOperator{\Quot}{Quot}

\DeclareMathOperator{\res}{{\rm res}}

\DeclareMathOperator{\SL}{{\rm SL}}
\DeclareMathOperator{\lexleq}{{\trianglelefteq}}
\DeclareMathOperator{\lexless}{{\vartriangleleft}}

\DeclareMathOperator{\dashV}{{\dashv\!\!\!\dashv}}
\DeclareMathOperator{\dotcup}{\dot{\cup}}
\newcommand{\lf}{\lfloor}
\newcommand{\rf}{\rfloor}
\newcommand{\blt}{\bullet}
\newcommand{\ovl}{\overline}

\newcommand{\wtilde}{\widetilde}
\newcommand{\ssstyle}{\scriptscriptstyle}
\newcommand{\talpha}{\widetilde{\alpha}}
\newcommand{\tbeta}{\widetilde{\beta}}
\newcommand{\BS}{\widehat\Sigma}
\newcommand{\surjto}{\mbox{$\hspace*{0.2em}\rightarrow\hspace*{-.8em}%
                                            \rightarrow\hspace*{0.2em}$}}
\newcommand{\injto}{\hookrightarrow}
\newcommand{\isoto}{\overset{\sim}{\longrightarrow}}
\newcommand{\tto}{\longrightarrow}
\newcommand{\Res}[2]{{#1}_{|_{\scriptstyle #2}}}
\newtheorem{lemma}{Lemma}[section]
\newtheorem{prop}[lemma]{Proposition}
\newtheorem{thm}[lemma]{Theorem}
\newtheorem{cor}[lemma]{Corollary}
\theoremstyle{remark}
\newtheorem*{rem}{Remark}
\newtheorem*{rems}{Remarks}
\theoremstyle{definition}

\numberwithin{equation}{section}

\title{The $T$--equivariant Cohomology of Bott--Samelson varieties}
\author{Martin H\"arterich}
\address{Bergische Universität Wuppertal\\
Fachbereich C, Mathematik und Naturwissenschaften\\
Gauß\-straße 20\\
42109 Wuppertal\\
Germany}
\curraddr{Dannheckerstraße 19\\ 69190 Walldorf\\ Germany}
\email{haerterich@math.uni-wuppertal.de}
\thanks{This work has been partially supported by the EC TMR network
``Algebraic Lie Representations'', contract no. ERB FMRX--CT97--0100.}
\date{\today}
\begin{document}
\maketitle
\begin{abstract}
We consider the $T$--equivariant cohomology of Bott--Samelson desingularisations
of Schubert varieties in the flag manifold of a connected semi--simple complex
algebraic group of adjoint type with maximal torus $T$.
We construct a combinatorially pure (in the sense of T.~Braden and R.~Macpherson)
sheaf on the Bruhat graph of the associated Weyl group such that its global
sections give the $T$--equivariant cohomology of our Bott--Samelson resolution.
\end{abstract}
%
%
\section*{Introduction}
%
In their paper \cite{BM01} Tom Braden and Robert Macpherson studied the notion
of sheaves on moment graphs and their relation to equivariant intersection cohomology
of certain complex projective varieties with torus actions.

They already pointed out that this approach is suited very well for studying the
$T$-equivariant intersection cohomology of Schubert varieties in the Flag variety
$G/B$ where $G$ is a connected semi--simple complex algebraic group of adjoint type,
$T$ is a maximal torus and $B$ is a Borel subgroup.
In this case the moment graph is the Bruhat graph of the Weyl group $W$,
i.~e.\ the oriented graph with vertices $W$ such that there is an arrow
from $x$ to $y$ whenever $y<x$ and $y=tx$ for some (not neccessarily simple)
reflection $t\in W$.

It turns out that some description of the $T$--equivariant cohomology ring
$\rH_T^\blt(G/B)$ known before (see e.~g.\ \cite{Br98}) may be regarded as the
ring of global sections of a 'constant' combinatorial sheaf on the Bruhat graph.

Using weight filtrations Braden and Macpherson were able to show much more.
They introduced combinatorially pure sheaves and proved that the $T$-equivariant
{\em intersection\/} cohomology of Schubert varieties in $G/B$ is obtained as the
global sections of indecomposable pure sheaves on the Bruhat graph

Now let $\pi\colon\BS\surjto \ovl{BwB/B} \injto G/B$ be a Bott--Samelson
desingularisation of some Schubert variety $\ovl{BwB/B}$.
(See the next section for more details.)
The motivation for the present work comes from the fact that as a consequence
of the decomposition theorem for equivariant perverse sheaves (see \cite{BL94}), $\rH_T^\blt(\BS)$ is the direct sum of the $T$--equivariant intersection cohomology
of some Schubert varieties.
Hence it equals the global section of the direct sum of certain combinatorially pure
sheaves corresponding to these Schubert varieties.

In this article we give a direct approach to show that $\rH_T^\blt(\BS)$ is the global sections of a pure sheaf on the associated Bruhat graph.
To achieve this we combine the combinatorics of galleries with localization techniques
for equivariant cohomology using results of A. Arabia on equivariant Euler classes
(see \cite{Ar98}) .

We give an overview of the organization of the article and the main results:

Section~\ref{Sec:Galleries} is a summary of basic concepts about Bott--Samelson
varieties and combinatorial galleries.
Section~\ref{Sec:Fibre} gives a description of the fibre $\BS_y$ of the
Bott--Samelson resolution over some $T$--fixed point $yB/B$ of the Schubert variety.
In particular, in Proposition~\ref{Prop:Fibre} we describe its decomposition into cells inherited from the Bialynicki--Birula decomposition of $\BS$.

In Section~\ref{Sec:Localization} we apply A. Arabia's equivariant Euler classes to our situation and describe some of them in terms of galleries. Furthermore we explore the
relations between the equivariant cohomology of $\BS$ and that of $\BS_x$.
Section~\ref{Sec:FixedPoints} describes the fixed point locus of the action of some subtori
of $T$ of codimension one and identifies them with Bott--Samelson varieties
for $G=\SL_2$, i.~e.\ with a product of projective lines.

This case ($G=\SL_2$) is treated in more detail in Section~\ref{Sec:SL2Case}.
We give two rather explicit descriptions of $\rH_T(\BS)$ here using the operation
of folding galleries.
We use these descriptions in Section~\ref{Sec:GeneralCase} to obtain results for the
general case.
In Theorem~\ref{Thm:HTBS1} we describe the image of the restriction to fixed points
$\rH_T(\BS)\injto\rH_T(\BS^T)$ a set of by non--linear congruences involving some
statistics on galleries. 
In Theorem~\ref{Thm:HTBS2} we give a description related to foldings.

Section~\ref{Sec:Sheaves} is devoted to an introduction into Braden and Macphersons
pure sheaves and a reinterpretation of the results of Section~\ref{Sec:GeneralCase}
in these terms.
In fact, we see that we have constructed a sheaf on the Bruhat graph of $W$ whose
module of global sections is just $\rH_T(\BS)$ (Theorem~\ref{Thm:BSSheaf}) and,
furthermore, that this sheaf is pure (Theorem~\ref{Thm:BSPurity}).
Finally, in Section~\ref{Sec:CombBases} we demonstrate how to calculate recursively
bases for the equivariant cohomology of the fibres $\Gamma_x$ and thus sheaves and
their decompositions.

The author would like to thank St\'ephane Gaussent for lots of discussions on the
subject and for even more motivation.
I would also like to thank Markus Reineke very much.
%
\section{Combinatorial Galleries}
\label{Sec:Galleries}
We fix a pinning (cf.~\cite{Ti87}) of our group $G$.
In particular $T$ is a fixed maximal torus and $R$ resp.\ $\Delta$ is the set of all
roots (resp.\ of the simple roots).
Denote by $B$ the Borel subgroup corresponding to our choice of simple
roots $\Delta$ and by $W$ the Weyl group generated by the simple reflections $S$.
Let $\leq$ be the Bruhat order and $\ell$ be the length function on $W$. 

For each root we have the one--parameter subgroup $p_\alpha\colon\C\to U_\alpha$.
Their mutual commutation laws can be found in \cite{Ti87}. Since we will need them only
very rarely, we cite them when we use them (see equations \eqref{Eq:ComRel1} and
\eqref{Eq:ComRel2}).

Denote by $P_\alpha$ the parabolic generated by $B$ and $s_\alpha$ and
(for any root $\alpha$) by $G_\alpha$ the subgroup generated by $T$ and $U_{\pm\alpha}$.
Set $B_\alpha=B\cap G_\alpha$.

We fix a sequence ${\bf s}=(s_1,\ldots,s_r)$ of simple reflections and set $w=s_1\cdots s_r$.
For $i=1,\ldots,r$ let $\alpha_i$ be the (simple) root corresponding to $s_i$.
Then the Bott--Samelson variety corresponding to ${\bf s}$ is
\[
    \BS = \BS({\bf s})
        = P_{\alpha_1}\times_B P_{\alpha_2}\times_B \cdots \times_B P_{\alpha_r}/B   \:,
\]
i.~e.\ the quotient of
$P_{\alpha_1}\times P_{\alpha_2}\times \cdots \times P_{\alpha_r}$ by the right
action of $B\times B\times \cdots \times B$ given by
$(p_1,p_2,\ldots,p_r)\cdot(b_1,b_2,\ldots,b_r)
=(p_1 b_1,b_1^{-1}p_2b_2,\ldots,b_{r-1}^{-1}p_rb_r)$.
It is a smooth variety of dimension~$r$.
We use the notation $[p_1,p_2,\ldots,p_r]$ for the point corresponding to the class of
$(p_1,p_2,\ldots,p_r)$.
Let $\pi\colon\BS\surjto G/B$ be the 'multiplication' map
$[p_1,p_2,\ldots,p_r]\mapsto p_1\cdots p_r B/B$ and denote by
$\BS_y=\pi^{-1}(yB/B)$ its fibre over the point $yB/B$.
If $s_1\cdots s_r$ is a reduced expression for $w$ then $\pi$ is a desingularisation map
onto its image which in this case is precisely the Schubert variety $\ovl{BwB/B}$.
\begin{rem}
The Bott--Samelson variety may also be seen as a variety of galleries inside the building
associated to the group $G$ (see \cite{CC83} or \cite{Ga01} for details).
\end{rem}
For $i=1,\ldots,r$ let $\gamma^i=\gamma_1\gamma_2\cdots\gamma_i$ and let
$\beta_i=\gamma^i(-\alpha_i)$.
The set of $T$--fixed points in $\BS$ is
\[
    \Gamma = \Gamma_{\bf s}
           = \{[\gamma_1,\ldots,\gamma_r]\mid \gamma_i=s_i \text{ or } \id \}
\]
Elements of this set are called {\em (combinatorial) galleries.\/}
They may be depicted as paths traversing the Weyl chambers (in the order
$\gamma^0=\id$, $\gamma^1$, $\gamma^2$,\ldots,$\gamma^r$) in such a way that they have
exactly one point in common with the wall corresponding $\beta_i$ when they move from
$\gamma^{i-1}$ to $\gamma^i$. 

Let $\Gamma_y=\Gamma\cap\BS_y$ be the set of galleries ending in $y$.
$\gamma_i$ is called a {\em bend\/} resp.\ a {\em crossing\/} if it equals $\id$ resp.\ $s_i$.
Then $\beta_i=\gamma^i(-\alpha_i)$ is the root corresponding to the wall of the $i$--th
crossing resp.\ bend of $\gamma$. Set $M_\alpha(\gamma)=\{i\mid\beta_i=\pm\alpha\}$.
Then the wall corresponding to $\beta_i$ is called load--bearing at $i$ if $\beta_i$
is positive. By abuse of language we will then also call $i$ or even $\gamma_i$
load--bearing. In terms of paths this is the case exactly if $\gamma$ is crossing
or bending {\em away\/} from the fundamental chamber.
Let $J(\gamma)$ be the set of load--bearing indices for $\gamma$ and
$J_\alpha(\gamma)=\{i \mid \beta_i=\alpha\}=J(\gamma)\cap M_\alpha(\gamma)$ the
load--bearing indices at the $\alpha$--wall.
Note that the assignment $\gamma\mapsto J(\gamma)$ is a bijection of $\Gamma$
onto the power set of $\{1,\ldots,r\}$.

Set $\wtilde{\alpha}_i=\gamma_i(-\alpha_i)$ and
$\wtilde{\beta}_i=\gamma^i(\wtilde{\alpha}_i)=\gamma^{i-1}(-\alpha_i)$.
Then $\gamma_i$ (resp.\ $i$) is called {\em defect\/}
if $\wtilde{\beta}_i$ is positive (equivalently, if $\gamma^{i-1}>\gamma^{i-1}s_{\alpha_i}$).
Let $D(\gamma)$ be the set of defect indices for $\gamma$ and
$D_\alpha(\gamma)=\{i \mid \wtilde{\beta}_i=\alpha\}=D(\gamma)\cap M_\alpha(\gamma)$ the
defect indices at the $\alpha$--wall.

Let $\sim_\alpha$ be the equivalence relation defined by
\[
 \delta\sim_\alpha\gamma :\!\iff \delta_i=\gamma_i \text{ unless } \beta_i=\pm\alpha  \:.
\]
\begin{rems}
(1)
$M_\alpha(\gamma)$ does not depend on $\gamma$ itself but only on its
$\sim_\alpha$--equivalence class.
\\
(2)
The $\alpha_i$ depend only on $\bf s$ whereas the $\talpha_i$, $\beta_i$ and $\tbeta_i$
depend also on the gallery~$\gamma$.
\end{rems}
\begin{lemma}
(a)
Let $i_1<i_2<\ldots<i_\ell$ be the elements of $M_\alpha(\gamma)$.
Then for $1\leq j<\ell$ we have
\[
      i_j\in J_\alpha(\gamma) \iff i_{j+1}\in D_\alpha(\gamma)
\]
and $i_\ell\in J_\alpha(\gamma)$ if and only if $\pi(\gamma)>s_\alpha\pi(\gamma)$.
\\
In particular, if $\gamma\sim_\alpha\delta$ and $\pi(\gamma)=\pi(\delta)$, then
$J_\alpha(\gamma)\subseteq J_\alpha(\delta) \iff D_\alpha(\gamma)\subseteq D_\alpha(\delta)$.
\\
(b)
$\# J(\gamma)- \# D(\gamma)=\ell(\pi(\gamma))$.
\end{lemma}
\begin{proof}
(a)
is left to the reader.
\\
(b)
By part (a) $\# J(\gamma)- \# D(\gamma)$ equals the number of positive roots $\alpha$ such that
$\pi(\gamma)>s_\alpha\pi(\gamma)$.
\end{proof}
We view $\BS$ as a closed subvariety of $(G/B)^r$ via the $T$--equivariant embedding
\begin{align*}
   \iota\colon\BS &\injto (G/B)^r \\
      [p_1,p_2,\ldots,p_r] &\mapsto (p_1,p_1p_2,\ldots, p_1p_2\cdots p_r) \:.
\end{align*}
The image of $\iota$ consists of all $(g_1,g_2,\ldots,g_r)$ such that
$(g_{i-1})^{-1}g_i \in P_{\alpha_i}$ for $i=1,2,\ldots,r$ (where $g_0=\id$).
If $G=\SL_2$ then $\iota$ is an isomorphism.

The Bott--Samelson variety has an open covering by sets
\[
     U^\gamma
   = \{[p_{\wtilde{\alpha}_1}(x_1)\gamma_1,\ldots,p_{\wtilde{\alpha}_r}(x_r)\gamma_r]\} \:.
\]
If $\alpha$ is a simple root, we have the relation
\begin{equation}\label{Eq:ComRel1}
  s_\alpha p_\beta(x) s_\alpha^{-1} = p_{s_\alpha(\beta)}(\pm x)
\end{equation}
the sign depending on the pinning of $G$ (see \cite{Ti87}).
Therefore we have
\begin{align}\label{Eq:iotap}
     &\iota([p_{\wtilde{\alpha}_1}(x_1)\gamma_1,...,p_{\wtilde{\alpha}_r}(x_r)\gamma_r]) \\
     &\qquad =\bigl(p_{\beta_1}(x_1)\gamma^1,p_{\beta_1}(x_1)p_{\beta_2}(\pm x_2)\gamma^2,
               \ldots,p_{\beta_1}(x_1)\cdots p_{\beta_r}(\pm x_r)\gamma^r\bigr) \:. \notag
\end{align}

Take a one parameter subgroup $\chi$ of $T$ such that the natural pairing of $\chi$ with
any $\alpha\in R^+$ is positive and consider the action of $\C^*$ on $\BS$ defined via
$\chi$.
As a consequence of \eqref{Eq:iotap} the Bialynicki--Birula decomposition of $\BS$ for
this action has cells
\begin{equation}\label{Def:Cgamma}
   C^\gamma=\{[p_{\wtilde{\alpha}_1}(x_1)\gamma_1,...,p_{\wtilde{\alpha}_r}(x_r)\gamma_r]
                    \mid x_i=0 \text{ if } \beta_i<0\}
\end{equation}
Note that this equals
$\{[p_{\wtilde{\alpha}_1}(x_1)\gamma_1,...,p_{\wtilde{\alpha}_r}(x_r)\gamma_r]
\mid x_i=0 \text{ unless } i\in J(\gamma)\}$.

For later use, we define some relations on the set of all combinatorial galleries:
\begin{align*}
 \delta\lexless\gamma
    &:\!\iff \exists i_0\colon \delta^{i_0}<\gamma^{i_0}
             \text{ and } \forall  i<i_0\colon \delta^i=\gamma^i      \\
 \delta\lexless_\alpha\gamma
    &:\!\iff \delta\lexless\gamma \text{ and } \delta\sim_\alpha\gamma \\
 \delta < \gamma
    &:\!\iff \exists i_0\colon \delta^{i_0}<\gamma^{i_0}
             \text{ and } \forall  i>i_0\colon \delta^i=\gamma^i      \\
 \delta\dashv\gamma
    &:\!\iff \delta\in\ovl{C_\gamma} \\
 \delta\dashV\gamma
    &:\!\iff \exists \delta_0,\delta_1,\ldots\colon\:
             \delta=\delta_0\dashv\delta_1\dashv\cdots\dashv\gamma
\end{align*}
As usual let us write $\delta\lexleq\gamma$ if $\delta\lexless\gamma$ or $\delta=\gamma$ etc.

\begin{rems}
(1) By considering the composition of the embedding $\iota$ with the projections onto the
various factors $G/B$ of $(G/B)^r$ one deduces that
$\delta\dashv\gamma$ (and by induction also $\delta\dashV\gamma$)
implies $\delta^i\leq\gamma^i$ for $i=1,2,\ldots,r$.
\\
(2) $\dashV$, $\lexleq$ and $<$ are orders. To prove it for $\dashV$ Remark (1) is useful.
\\
(3)
Although the definitions of $\lexless$ and $<$ look very similar the order $\lexleq$ is
total, but $<$ is not.
Only its restriction to $\Gamma_x$ (for some $x\in W$) is a total order.
(We refer to $\lexleq$ as the {\em lexicographic order\/} and use it mainly to order
galleries when used as indices of matrices.)
\\
(4) $\dashv$ is not transitive in general.
\end{rems}
The following implications are clear from the definitions and the remarks:
\begin{align}
  &\delta\dashv\gamma \Longrightarrow \delta\leq\gamma     \label{Eq:ClosureBruhat} \\
  &\delta\dashv\gamma \Longrightarrow \delta\lexleq\gamma  \label{Eq:ClosureLex}    \\
  &J(\delta)\subseteq J(\gamma) \Longrightarrow \delta\lexleq\gamma \\
  &\delta\sim_\alpha\gamma \text{ and } J_\alpha(\delta)\subseteq J_\alpha(\gamma) 
                 \Longrightarrow \delta\lexleq_\alpha\gamma \\
  &\delta\lexleq_\alpha\gamma \Longrightarrow \delta\lexleq\gamma
\end{align}
%
\section{Description of the fibre}
\label{Sec:Fibre}
Fix some $T$--fixed point in $G/B$ i.~e.\ an element $y$ of the Weyl group $W$.
The fibre $\BS_y$ inherits a cell decomposition by intersecting with the
Bialynicki--Birula cells of $\BS$.
Clearly, $\BS_y\cap C^\gamma = \emptyset$ unless $\pi(\gamma)=y$.
For $\gamma\in\Gamma_y$ let $C^\gamma_y=\BS_y\cap C^\gamma$.

\begin{prop}\label{Prop:Fibre}
$C^\gamma_y$ is determined as a subvariety of $C^\gamma$ by $\ell(y)$ equations
$X_\alpha=0$, $\alpha\in R^+\cap y(R^-)$,
where $X_\alpha=\Bigl(\sum\limits_{i\colon \beta_i=\alpha} \pm\; x_i \Bigr) + \wtilde{X}_\alpha$
for a polynomial $\wtilde{X}_\alpha$ of degree at least $2$ in the variables
$\{x_j\mid\Ht(\beta_j)<\Ht(\alpha))\}$.
(The signs are determined by the pinning of $G$.)
\\
In particular, $C^\gamma_y$ is the graph of a polynomial function in $\#J(\gamma)-\ell(y)$ variables
and hence a cell of dimension $\#D(\gamma)$.
\end{prop}
\begin{proof}
Let $p = [p_{\wtilde{\alpha}_1}(x_1)\gamma_1,\ldots,p_{\wtilde{\alpha}_r}(x_r)\gamma_r]$
be an element of $U^\gamma$.
We already know that
$\pi(p) = p_{\beta_1}(x_1)p_{\beta_2}(\pm x_2)\cdots p_{\beta_r}(\pm x_r)yB/B$.
Now assume that furthermore $p \in C^\gamma$, i.~e.\ that $x_i=0$ whenever the root
$\beta_i$ is negative.
We want to apply commutation relations among the $p_\alpha$ with $\alpha>0$ to get
an expression of the form
\[
    \pi(p)
  = \Bigl(\prod_{\alpha\in R^+\cap y(R^-)} p_{\alpha}(X_\alpha)\Bigr)
         \Bigl(\prod_{\alpha\in R^+\cap y(R^+)} p_{\alpha}(X_\alpha)\Bigr) yB \:.
\]
(Here it is understood that we have fixed some order of the factors $p_{\alpha}(X_\alpha)$
in both products.)
In fact we can achieve this by using the relations
\begin{equation}\label{Eq:ComRel2}
    p_{\alpha_1}(x_1)p_{\alpha_2}(x_2)
  = p_{\alpha_2}(x_2)\Bigl(\prod_{(m,n)\colon m\alpha_1+n\alpha_2\in R}p_{m\alpha_1+n\alpha_2}
      (C_{m,n}x_1^mx_2^n)\Bigr)p_{\alpha_1}(x_1)
\end{equation}
for certain integers $C_{m,n}$ depending on $m$, $n$ and also on $\alpha_1$, $\alpha_2$
(see \cite{Ti87}).

We conclude that $X_\alpha$ is a polynomial in the variables
$x_1$,\ldots, $x_r$ with integral coefficients.
Moreover, the variable $x_i$ is involved linearly only in the argument of $p_{\beta_i}$
and non--linearly only in the arguments of some $p_\alpha$ for $\Ht(\alpha)>\Ht(\beta_i)$.
Therefore,
\[
    \pi(p)
 = \Bigl(\prod_{\alpha\in R^+\cap y(R^-)} p_{\alpha}(X_\alpha)\Bigr) yB/B  \\
 = y \Bigl(\prod_{\alpha\in R^+\cap y(R^-)} p_{y^{-1}(\alpha)}(\pm X_\alpha)\Bigr) B/B \:,
\]
where the $X_\alpha$ are as in the proposition.
This equals $yB/B$ if and only if all $X_\alpha$ vanish.

The last statement of the proposition is clear.
\end{proof}
%
\section{Localization and inverse Euler classes}
\label{Sec:Localization}
Let $\K$ be a field or $\K$ be the ring of integers $\Z$.
Let $\rH_T^\blt(\BS)$ be the $T$--equivariant cohomology with coefficients in $\K$
of the Bott--Samelson variety and $A=\rH_T^\blt(pt)$ the equivariant cohomology ring of
a point. Denote by $Q=\Quot A$ its quotient field.

By Chern--Weil theory $A$ is known to be the ring $\K\otimes_\Z\Sym\langle\Delta\rangle$
where $\Sym\langle\Delta\rangle$ is the symmetric algebra (over $\Z$) of the root lattice.

Since $\BS$ is smooth and projective there is a filtration by compact subsets
$\BS=X_n\supseteq X_{n-1} \supseteq \ldots \supseteq X_0=\emptyset$ such that
each $X_i\setminus X_{i-1}$ is a Bialynicki-Birula cell (and hence of even real
dimension).
Therefore, the (non--equivariant) cohomology $\rH^\blt(\BS)$ vanishes in odd degrees
and consequently the spectral sequence of equivariant cohomology degenerates at the
$E_2$--term, i.~e.\ we have 
\begin{equation} \label{Eq:SpecSeq}
    \rH_T^\blt(\BS) = A \otimes \rH^\blt(\BS)
\end{equation}
as a graded $A$--algebra (see \cite[Theorem 5.2]{McC01}).
In particular, $\rH_T^\blt(\BS)$ is a free $A$--module.

The Bialynicki-Birula cells $C^\gamma$ define a basis of the ordinary homology of $\BS$.
Let $\{\mu_\delta^0 \mid \delta\in\Gamma\}\subseteq\rH^\blt(\BS)$ be the corresponding
dual basis in the (non--equivariant) cohomology.
Now by \eqref{Eq:SpecSeq} we have an $A$--basis $\{1\otimes\mu_\delta^0\mid\delta\in\Gamma\}$
of $\rH_T^\blt(\BS)$ and as in \cite[2.5]{Ar89} we can apply a Gram--Schmidt algorithm 
to get a basis $\{\mu_\delta\mid\delta\in\Gamma\}\subseteq\rH_T^\blt(\BS)$ such that
\begin{equation}  \label{Eq:Integration}
     \int_{\ovl{C^\gamma}}\mu_\delta
   = \begin{cases}
        1, &\text{ if } \gamma=\delta, \\
        0, &\text{ else.}
     \end{cases}
\end{equation}
Here $\int_{\ovl{C^\gamma}}\colon\rH_T^\blt(\BS)\to A$ is the $A$--linear map (of degree
$-\dim_\R(C^\gamma)$) obtained by composing the restriction
$\rH_T^\blt(\BS)\to\rH_T^\blt(\ovl{C^\gamma})$ with the integration defined
in \cite[1.4]{Ar98}.

Let $\Omega_\Gamma$ be the set of all maps $f\colon\Gamma\to A$. This is a ring under
pointwise addition and multiplication of functions.
The inclusion of fixed points $i\colon\BS^T\injto\BS$ induces a restriction homomorphism
$i^*\colon\rH_T^\blt(\BS)\tto\rH_T^\blt(\BS^T)$.
It is injective and becomes an isomorphism when tensored with $Q$.
Our aim is to describe the image of this homomorphism.
Because $\rH_T^\blt(\BS^T)=\bigoplus_{\gamma\in\Gamma}\rH_T^\blt(\gamma)$ can
canonically be identified with $\Omega_\Gamma$ we will view $\rH_T^\blt(\BS)$
as a subring of $\Omega_\Gamma$.

Let $V=\rH_T^\blt(\BS^T)\otimes_A Q=\bigoplus_{\gamma\in\Gamma} V_\gamma$ where
$V_\gamma=\rH_T^\blt(\gamma)\otimes_A Q\cong Q$ (as a $Q$--algebra).
Under the natural identification the projections $\pi_\gamma\colon V\surjto V_\gamma$
form a basis of $V^*$ dual to the basis $(1_\gamma\otimes 1)_{\gamma\in\Gamma}$ of $V$.
(Here $1_\gamma$ is the $1\in \rH_T^\blt(\gamma)=A$.)
\begin{thm}[Arabia] \label{Thm:Arabia}
(a) There are certain (T-equivariant) Euler classes $\Eu_T(\delta,\ovl{C^\gamma})\in Q$
such that denoting
\begin{equation} \label{Eq:Integral}
    I_\gamma=\sum_{\delta\dashv\gamma}\frac{\pi_\gamma}{\Eu_T(\delta,\ovl{C^\gamma})}\in V^*
\end{equation}
we have $\int_{\ovl{C^\gamma}}=I_\gamma\circ i^*$.
\\
(b) $\Eu_T(\delta,\ovl{C^\gamma})$ depends only on a $T$--invariant neighbourhood of the fixed point $\delta$.
\\
(c) If $\delta\dashv\gamma$ is a smooth point of $\ovl{C^\gamma}$ then $\Eu_T(\delta,\ovl{C^\gamma})$
is the Pfaffian determinant of the induced action of $T$ on the tangent space $T_\delta\ovl{C^\gamma}$.
\end{thm}
\begin{proof}[Remarks on the proof]
(a) is the localization formula of \cite[2.3.1]{Ar98}. The fact that the Euler classes are non--zero
is proved in \cite[2.2.1]{Ar98}.
(b) is \cite[1.4.1]{Ar98}.
(c) is \cite[2.4.1]{Ar98}.
\end{proof}
Let $\gmax$ be the gallery corresponding to the cell of maximal dimension (i.~e.\ all walls
of $\gmax$ are load--bearing). Then $\ovl{C^{\gmax}}=\BS$ and $U^\gamma$ is a neighbourhood
of $\gamma$ in $\ovl{C^{\gmax}}$ isomorphic to $\C^r$.
Similarly, $C^\gamma$ is a neighbourhood of $\gamma$ in $\ovl{C^\gamma}$ isomorphic to $\C^{\#J(\gamma)}$.
Using part (c) of the theorem we obtain the following explicit formulas:
\begin{align}
      \Eu_T(\gamma,\ovl{C^{\gmax}})
   &= \prod_{i\in J(\gmax)} \gamma^i(-\alpha_i)
    = \prod_{i=1}^r \beta_i
    = \prod_{\alpha\in R^+} (-1)^{\#J_\alpha(\gamma)}\cdot(-\alpha)^{\#M_\alpha(\gamma)}   \label{Eq:EuTg0}
\\
      \Eu_T(\gamma,\ovl{C^\gamma})
   &= \prod_{i\in J(\gamma)}   \gamma^i(-\alpha_i)
    = \prod_{i\colon \beta_i>0} \beta_i
    = \prod_{\alpha\in R^+} \alpha^{\#J_\alpha(\gamma)}                \label{Eq:EuTg}
\end{align}
Let
\[
      E=\left(\frac{1}{\Eu_T(\delta,\ovl{C^\gamma})}\right)_{\gamma,\delta\in\Gamma}
\]
be the matrix of inverse Euler classes, where it is understood that $ \frac{1}{\Eu_T(\delta,\ovl{C^\gamma})}=0$
unless $\delta\dashv\gamma$.
Let
\[
     H=\bigl(\Res{\mu_\gamma}{\delta}\bigr)_{\delta,\gamma\in\Gamma}
\]
be the matrix with columns the restrictions of the basis vectors $\mu_\gamma$ to the $T$--fixed points.
\begin{cor}
$H=E^{-1}$, in particular the diagonal elements are
$\Res{\mu_\gamma}{\gamma}=\prod_{i\colon\beta_i>0}\beta_i
=\prod_{\alpha\in R^+} \alpha^{\#J_\alpha(\gamma)}$.
\end{cor}
\begin{proof}
We may rephrase part (a) of the theorem as follows: $(I_\gamma)_{\gamma \in\Gamma}\subset V^*$
is the basis dual to the basis $(i^*\mu_\delta\otimes 1)_{\gamma \in\Gamma}\subset V$.
Note that by \eqref{Eq:Integral} the transpose ${}^t\!E$ is just the base
change matrix between the two bases $(I_\gamma)_{\gamma\in\Gamma}$ and
$(\pi_\delta)_{\delta\in\Gamma}$ of $V^*$.
Consequently, $E$ expresses the base change between their respective dual bases 
$(i^*\mu_\gamma)_{\gamma\in\Gamma}$ and $(1_\delta\otimes 1)_{\delta\in\Gamma}$ of $V$.

Since $E$ is a lower triangular matrix, the diagonal elements of $H$ can be calculated
using~\eqref{Eq:EuTg}.
\end{proof}
\begin{rem}
Unfortunately, it is hard to give a combinatorial description of the relation
$\delta\dashv\gamma$ in general and to find the corresponding inverse Euler classes
at singular points of $\ovl{C^\gamma}$.
Therefore we will use this direct approach only in the case of $G=\SL_2$.
\end{rem}

Note that the preceeding arguments are true for any order chosen for the indices of the
matrices $E$ and $H$.
If we take the lexicographic order $\lexleq$ then by \eqref{Eq:ClosureLex} $E$
and hence $H$ are lower trigonal.
By \eqref{Eq:ClosureBruhat} the same is true for any (total) order compatible with $\le$.
For the rest of this section let us fix such an order.

%
%

From the discussion in Section~\ref{Sec:Fibre} we know that (even if  $\BS_y$ is not
smooth) there is a filtration by compact subsets
$\BS_y=X_n\supseteq X_{n-1} \supseteq \ldots \supseteq X_0=\emptyset$ such that
each $X_i\setminus X_{i-1}$ is a cell of even real dimension.
Therefore, we may repeat all the arguments for the fibres $\BS_y$ instead of $\BS$:
For each $y\in W$ the $T$--equivariant cohomology $\rH_T^\blt(\BS_y)$ of the fibre
over $y$ is a free $A$--module and a basis $(\mu_{\gamma,y})_{\gamma\in\Gamma_y}$
with a property analogous to \eqref{Eq:Integration} is obtained from the
Bialynicki--Birula decomposition of $\BS_y$.

For $y\in W$ let $d_y$ be the product of the weights of the $T$--action on the Bruhat cell $ByB/B$,
i.~e.\ $d_y=\prod\limits_{\alpha\colon s_\alpha(y)<y}\alpha$.
Since $y$ is a smooth point of the Schubert variety $\ovl{ByB/B}$ we have
\[
  \Eu_T(y,\ovl{ByB/B}) = \Eu_T(y,ByB/B) = d_y \:.
\]
Consider $\gamma,\delta\in\pi^{-1}(y)$ such that $\delta\dashv\gamma$.
Then $\ovl{C^\gamma}\cap\pi^{-1}(ByB/B)$ is a $T$--invariant neighbourhood of $\delta$ in $\ovl{C^\gamma}$
isomorphic to $\ovl{C^\gamma_y} \times ByB/B$.
Therefore we can apply \cite[2.6.1--1]{Ar98} to see that
\begin{equation} \label{Eq:EuTProd}
    \Eu_T(\delta,\ovl{C^\gamma}) = \Eu_T(\delta,\ovl{C^\gamma_y})\cdot\Eu_T(y,ByB/B)
                                 = d_y\cdot\Eu_T(\delta,\ovl{C^\gamma_y})   \:.
\end{equation}
Let
\begin{equation} \label{Eq:DefHy}
     H_y=\bigl(\Res{\mu_{\gamma,y}}{\delta}\bigr)_{\delta,\gamma\in\Gamma_y}
\end{equation}
be the matrix with columns the restrictions of the basis vectors $\mu_{\gamma,y}$
to the $T$--fixed points in the fibre $\BS_y$.
\begin{prop} \label{Prop:EuTRes}
Let $y\in W$ and $\gamma$, $\delta\in\Gamma_y$.
Then $\Res{\mu_\gamma}{\delta} = d_y\cdot\Res{\mu_{\gamma,y}}{\delta}$,
in particular
\[
\Res{\mu_{\gamma,y}}{\gamma} = \prod\limits_{i\in D(\gamma)}\wtilde{\beta}_i
=\prod\limits_{\alpha\in R^+} \alpha^{\#D_\alpha(\gamma)}  \:.
\]
\end{prop}
\begin{proof}
For any $y\in W$ let $\wtilde{H}_y$ be the minor
$(\Res{\mu_\gamma}{\delta})_{\gamma,\delta\in\Gamma_y}$ of $H$.
We have to show that $\wtilde{H}_y=d_y \cdot H_y$.
This can be seen easily if we order the indices of the matrix in some way that is
compatible with the order $\leq$ on galleries.
Then the matrix $E$ is block diagonal where the diagonal blocks correspond to the
$\Bigl(\dfrac{1}{\Eu_T(\delta,\ovl{C^\gamma_y})}\Bigr)_{\gamma,\delta\in\Gamma_y}
=\wtilde{H}_y^{-1}$.
\end{proof}

%
%

Let $\gmax$ be the gallery corresponding to the cell of maximal dimension, i.~e.\ the one
with $J(\gmax)=\{1,2,\ldots,r\}$.
Then there is a non--degenerate pairing (Poincaré duality)
\begin{align*}
   \langle\:,\rangle_T\colon \rH_T^\blt(\BS) \times \rH_T^\blt(\BS) &\to A \\
      (\mu,\nu)                           &\mapsto \int_{\ovl{C^{\gmax}}} \mu \cup \nu   \:.
\end{align*}
Remember that we regard $\rH_T^\blt(\BS) \times \rH_T^\blt(\BS)$ as a submodule of $\Omega_Q\times\Omega_Q$ and let $D$ be the diagonal matrix with entries $\Eu_T(\delta,\ovl{C^\gmax})^{-1}$ ($\delta\in\Gamma $).
Then the integration formula in Theorem~\ref{Thm:Arabia}(a) tells us that
the Poincar\'e duality pairing is just the restriction of the pairing
$\Omega_Q\times\Omega_Q \to A$ with matrix $D$ with respect to the basis of $\Omega_Q$
consisting of the characteristic functions of the galleries in $\Gamma$.

Using \cite[1.4--1(c)]{Ar98} it follows that under the identification of \eqref{Eq:SpecSeq}
we have
\begin{equation} \label{Eq:PDPairing}
     \langle f\otimes\mu,g\otimes\nu\rangle_T = fg \cdot \langle\mu,\nu\rangle
\end{equation}
where $\langle\:,\rangle\colon \rH^\blt(\BS) \times \rH^\blt(\BS) \to \K$ is the usual
(non--equivariant) Poincar\'e duality pairing.
\begin{prop} \label{Prop:PDBasis}
For any $A$--basis of\/ $\rH_T^\blt(\BS)$ there is a dual basis with respect to the Poincar\'e duality pairing.
Specifically, the basis $\{\mu_\gamma^*\mid\gamma\in\Gamma\}$ dual to $\{\mu_\gamma\mid\gamma\in\Gamma\}$ is given by
\begin{equation} \label{Eq:PDBasis}
      \Res{\mu_\gamma^*}{\delta}=\frac{\Eu_T(\delta,\ovl{C^{\gmax}})}{\Eu_T(\delta,\ovl{C^{\gamma}})} \:,
\end{equation}
in particular: $\Res{\mu_\gamma^*}{\gamma}=\prod\limits_{i\colon\beta_i<0}\beta_i
=\prod\limits_{\alpha\in R^+}\alpha^{\#M_\alpha(\gamma)-\#J_\alpha(\gamma)}$.
\end{prop}
\begin{rem}
We implicitly claim that the expressions in the proposition are well-defined elements of $A$.
Since $\ovl{C^{\gmax}}=\BS$ is smooth this may also be seen from \cite[4.2]{Br97} or
\cite[Proposition 2.6--1(b)]{Ar98}.
\end{rem}
\begin{proof}
By standard linear algebra, it is clear that if there is {\em one\/} $A$--basis of 
$\rH_T^\blt(\BS)$ possessing a dual basis then the same is true for {\em any\/} $A$--basis.
But then we see from equation~\eqref{Eq:PDPairing} that the existence of a basis
$\{(\mu_\gamma^0)^*\}\subseteq\rH^\blt(\BS)$ dual to the basis $\{\mu_\gamma^0\}$
(i.~e.\ the non--equivariant Poincar\'e duality) implies the existence of a basis
$\{1\otimes{\mu_\gamma^0}^*\}\subseteq\rH_T^\blt(\BS)$ dual to the basis
$\{1\otimes\mu_\gamma^0\}$.

Now let
\begin{equation} \label{Eq:Hstar}
     H^* = \bigl(\Res{\mu_\gamma^*}{\delta}\bigr)_{\delta,\gamma\in\Gamma}
\end{equation}
and let $D$ be as above. Then ${}^tH^*\cdot D\cdot H$ is the identity matrix.
Thus $H^* = {}^t(H\cdot D)^{-1}= D^{-1} {}^t E$.
This gives the explicit formula.
For the ``in particular'' statement use \eqref{Eq:EuTg0} and \eqref{Eq:EuTg}.
\end{proof}
\begin{rem}
Note also that
\[
     H^{-1}H^*=E\cdot D^{-1}\cdot {}^t E
\]
is the matrix (obviously symmetric) that expresses the base change between
$\{\mu_\gamma^*\mid\gamma\in\Gamma\}$ and $\{\mu_\gamma\mid\gamma\in\Gamma\}$.
Hence it and its inverse both have entries in $A$.
(This holds even for integral coefficients.)
Moreover, the same argument works for {\em any\/} matrix $H$ such that the columns
represent a basis for $\rH_T^\blt(\BS)$.
\end{rem}
%
\section{Fixed points of Subtori of Codimension One}
\label{Sec:FixedPoints}
In this section we study the fixed point set of $\BS$ under the action of a codimension
one subtorus $T_\alpha=\ker(\alpha)\subset T$, where $\alpha$ is a fixed positive root.
We start with the following observation:
\begin{lemma} \label{L:CgFP}
The $T_\alpha$--fixed points in $C^\gamma$ are
\[
   (C^\gamma)^{T_\alpha}=\{[p_{\talpha_1}(x_1)\gamma_1,...,p_{\talpha_r}(x_r)\gamma_r]
                           \mid x_i=0 \text{ if } \beta_i\not=\alpha\} \:.
\] 
\end{lemma}
%
%
Denote by $\gmin$ the unique element in the $\sim_\alpha$--equivalence class of $\gamma$
that has no load--bearing $\alpha$--wall.
Let $i_1<i_2<\ldots<i_\ell$ be the elements of $M_\alpha(\gamma)$.
For $i=1,2,\ldots,r$ let $\lf i\rf^\alpha_\gamma$ be the largest element of
$M_\alpha(\gamma)\cup\{0\}$ that does not exceed $i$.
\begin{rem}
$\gamma^i=s_\alpha\gmin^i$ if $\lf i\rf^\alpha_\gamma\in J_\alpha(\gamma)$
and $\gamma^i=\gmin^i$ else.
\end{rem}

Now we want to define maps
$v_\gamma^\alpha \colon \bigl(G_\alpha/B_\alpha\bigr)^{\#M_\alpha(\gamma)} \to \BS$
such that $\iota \circ v_\gamma^\alpha$ is given by
\begin{equation} \label{Eq:assign}
          \bigl(G_\alpha/B_\alpha\bigr)^{\#M_\alpha(\gamma)} \to (G/B)^r,
                        \quad (g_1,g_2,\ldots) \mapsto (h_1,h_2,\ldots)  \:,
\end{equation}
where $h_i=g_j\,{\gmin^i}$ if $\lf i\rf_\gamma^\alpha=i_j$. (Here: $g_0=\id$)
\begin{rem}
Clearly, this does not depend on $\gamma$ itself but only on its
$\sim_\alpha$--equivalence class.
\end{rem}
To show that \eqref{Eq:assign} is well defined we have to check that
\begin{equation} \label{Eq:Commut}
     U_\alpha \gmin^i \subseteq \gmin^i B
\end{equation}
or equivalently $(\gmin^i)^{-1}(\alpha)\in R^+$ for all $i$.
But this is obvious for $i=0$ (here $\gmin^0=\id$) and follows
inductively for $i=1,2,\ldots$ because $(\gmin^i)^{-1}(\alpha)\in R^-$ {\em and\/} $(\gmin^{i-1})^{-1}(\alpha)\in R^+$ can occur only if $(\gmin)_i$ is a reflection
(hence $s_{\alpha_i}$) and $(\gmin^{i-1})^{-1}(\alpha)$ is the corresponding simple root
$\alpha_i$. But then $\gmin^i(-\alpha_i)=\gmin^{i-1}(\alpha_i)=\alpha$,
i.~e.\ $\gmin$ has a load--bearing wall. Contradiction!

Furthermore, the image of the assignment \eqref{Eq:assign} is in $\iota(\BS)$.
For this we have to show that $h_{i-1}^{-1}h_i\in P_{\alpha_i}$ for all $i$.
This is clear if $\lf i\rf_\gamma^\alpha=\lf i-1\rf_\gamma^\alpha$
Otherwise we have $i=i_j$ for some $j$. It follows that $\gmin^i(-\alpha_i)=-\alpha$
as well as $(\gmin)_i=\id$ and consequently
$h_{i-1}^{-1}h_i\in(\gmin^i)^{-1}G_\alpha\gmin^i \subseteq G_{(\gmin^i)^{-1}(\alpha)}
\subseteq P_{\alpha_i}$.

Finally, $v_\gamma^\alpha$ is injective.
For suppose that $\iota\circ v_\gamma^\alpha$ maps $(g_1,g_2,\ldots)$ and $(g'_1,g'_2,\ldots)$ to the same element of $(G/B)^r$.
Then $(g'_j \gmin^{i_j})^{-1}g_j\gmin^{i_j} \in B \cap G_{\alpha_{i_j}}$
for $j=1,2,\ldots$ (since $(g'_j)^{-1}g_j \in G_\alpha$ and $-\alpha=\gmin^{i_j}(-\alpha_{i_j}))$.
Therefore $(g'_j)^{-1}g_j\in \gamma^{i_j}B_{\alpha_{i_j}}(\gamma^{i_j})^{-1} = B_\alpha$.
\begin{prop} \label{Prop:vga2}
(a) $(C^\gamma)^{T_\alpha}=v_\gamma^\alpha(X_1\times X_2\times\cdots)$
where $X_j=U_\alpha s_\alpha B_\alpha/B_\alpha$ if $j\in J_\alpha(\gamma)$ and $X_j=B_\alpha/B_\alpha$
else.
\\
(b) $\im(v_\gamma^\alpha) \subseteq \iota(\BS^{T_\alpha})$
\\
(c) If $\gamma\not\sim_\alpha\delta$ then $\im(v_\gamma^\alpha) \cap \im(v_\delta^\alpha)=\emptyset$.
\end{prop}
\begin{proof}
(a)
By \eqref{Eq:iotap} and Lemma~\ref{L:CgFP} we have
\begin{align*}
    &  \iota\bigl((C^\gamma)^{T_\alpha}\bigr) \\
    &= \{(p_{\beta_1}(x_1)\gamma^1,p_{\beta_1}(x_1)p_{\beta_2}(x_2)\gamma^2,\ldots,
       p_{\beta_1}(x_1)\cdots p_{\beta_r}(x_r)\gamma^r)\mid x_i=0 \text{ if } i\not\in J(\gamma)\} \\
    &= \{(p_\alpha(x_1)\gamma^1,p_\alpha(x_1+x_2)\gamma^2,\ldots,
       p_\alpha(x_1+\ldots+x_r)\gamma^r)\mid x_i=0 \text{ if } \beta_i\not=\alpha\} \\
    &= \iota\bigl(v_\gamma^\alpha(X_1\times X_2\times\cdots)\bigr) \:.
\end{align*}
(For the last step use that $\gamma^i=\gmin^i$ if $\lf i\rf_\gamma^\alpha\not\in J_\alpha(\gamma)$
together with \eqref{Eq:Commut}.)
\\
(b) is clear by the above remarks.
\\
(c)
We may assume that $\gamma=\gmin$ and $\delta=\delta_{\rm min}$.
In particular, the smallest $i$ such that $\delta_i\not=\gamma_i$ does not correspond to an
$\alpha$--wall. Now suppose that $\im v_\gamma^\alpha \cap \im v_\delta^\alpha$ is non--empty.
This implies that there are $g$, $g'\in G_\alpha$ such that $g\gamma^iB=g'\delta^iB$.
This in turn leads to $(\delta^i)^{-1}(g')^{-1}g\gamma^i\in B$.
But $(\delta^i)^{-1}(g')^{-1}g\gamma^i\in G_{(\delta^i)^{-1}(\alpha)}s_{\alpha_i}$ which
does not contain any element of $B$ unless $(\delta^i)^{-1}(\alpha)=\pm\alpha_i$,
i.~e.\ unless $i$ corresponds to an $\alpha$--wall. 
\end{proof}
\begin{cor} \label{Cor:Closure}
$\ovl{(C^{\gamma})^{T_\alpha}}=v_\gamma^\alpha(Y_1\times Y_2\times\cdots)$
where $Y_j=G_\alpha/B_\alpha$ if $j\in J_\alpha(\gamma)$ and $Y_j=B_\alpha/B_\alpha$
else. Hence it is the disjoint union of all $(C^{\delta})^{T_\alpha}$
for all $\delta\sim_\alpha\gamma$ such that $J_\alpha(\delta)\subseteq J_\alpha(\gamma)$.
\end{cor}
\begin{proof}
Follows from part (a) of the proposition since $v_\gamma^\alpha$ is a proper map.
\end{proof}
\begin{cor}
The connected components of $\BS^{T_\alpha}$ are precisely all $\im(v_\gamma^\alpha)$
where $\gamma$ runs over a set of representatives for the $\sim_\alpha$--equivalence classes
of $\Gamma$. Hence they are all irreducible.
\end{cor}
\begin{proof}
By part (a) and (b) of the proposition
\[
             \BS^{T_\alpha}
   =         \bigcup_{\gamma\in\Gamma} (C^\gamma)^{T_\alpha}
   \subseteq \bigcup_{\gamma\in\Gamma} \im(v_\gamma^\alpha)
   \subseteq \BS^{T_\alpha} \:.
\]
Now the claim follows using (c) and the irreducibility of $(G_\alpha/B_\alpha)^{\#M_\alpha(\gamma)}$.
\end{proof}
%
\section{The $\SL_2$-Case}
\label{Sec:SL2Case}
As a first step towards the description of $\rH_T^\blt(\BS)$ we consider the case where $G=\SL_2$.
Let $\alpha$ be {\em the\/} simple root and $s$ be {\em the\/} simple reflection.
We have $G_\alpha=P_\alpha=G$ and $G/B$ is a projective line.

We have $(C^\gamma)^{T_\alpha}=C^\gamma$ because $T_\alpha=\{\pm 1\}$.
Furthermore $\iota\circ v_\gamma^\alpha$ is the identity map and $\iota$
is an isomorphism, in particular $\iota(\ovl{C^\gamma})=\ovl{\iota(C^\gamma)}$.
Hence we get from Corollary~\ref{Cor:Closure}
\begin{lemma}
$\iota(\ovl{C^\gamma})=Y_1\times\cdots\times Y_r$ where $Y_j=G/B$ if $j\in J(\gamma)$
and $Y_j=B/B$ else.
Therefore $\ovl{C^\gamma}$ is a smooth variety of dimension $\#J(\gamma)$
and $(\ovl{C^\gamma})^T =\{\delta\in\Gamma \mid J(\delta)\subseteq J(\gamma)\}$.
\end{lemma}
This enables us to write down the following explicit formulas.
\begin{prop} \label{Prop:EuT_HT_SL2}
(a)
Let $\delta$, $\gamma\in\Gamma$. Then
\[
      \Eu_T(\delta,\ovl{C^\gamma})
    = \begin{cases}
         (-1)^{\#J(\delta)}(-\alpha)^{\# J(\gamma)}\,,
                        &\text{if } J(\delta)\subseteq J(\gamma) \\
         \infty \,,     &\text{else.}
      \end{cases}
\]
(b) The basis $\{\mu_\gamma\mid\gamma\in\Gamma\}$ of $\rH_T^\blt(\BS)$ is determined by
\[
      \Res{\mu^\gamma}{\delta}
    = \begin{cases}
         \alpha^{\#J(\gamma)}, &\text{if } J(\delta)\supseteq J(\gamma) , \\
         0,                    &\text{else.}
      \end{cases}\:
\]
(c) The Poincar\'e dual basis $\{\mu_\gamma^*\mid\gamma\in\Gamma\}$ of $\rH_T^\blt(\BS)$
is given by
\[
      \Res{\mu_\gamma^*}{\delta}
    = \begin{cases}
         (-\alpha)^{r-\#J(\gamma)}, &\text{if } J(\delta)\subseteq J(\gamma) , \\
         0,                         &\text{else.}
      \end{cases}\:
\]
\end{prop}
\begin{rem}
Let us assume that the indices of $H$ and $H^*$ are ordered according to the
lexicographic order.
Then (b) and (c) show that the matrix $H^*$ (see~\eqref{Eq:Hstar}) is obtained
from $H$ by reversing the order of the rows and columns and substituting $(-\alpha)$
for $\alpha$.
\end{rem}
\begin{proof}
(a)
$U^\delta \cap \ovl{C^\gamma}$ is a $T$--invariant neighbourhood of $\delta$ in $\ovl{C^\gamma}$
isomorphic to $\C^{\#J(\gamma)}$. Therefore using Theorem~\ref{Thm:Arabia}(c) we obtain the formula
$\Eu_T(\delta,\ovl{C^\gamma})=\prod_{i\in J(\gamma)}\delta^i(-\alpha)$ for $J(\delta)\subseteq J(\gamma)$
whence the claim.
\\
(b)
It follows from (a) that the matrix of inverse Euler classes is
\[
   E = \begin{pmatrix} E' & 0 \\ -\frac{1}{\alpha}E' & \frac{1}{\alpha}E' \end{pmatrix}\:,
\]
where $E'$ is the matrix of inverse Euler classes for ${\bf s'}=(s,\ldots,s)$ (of length $r-1$)
instead of~${\bf s}$.
Therefore
\[
   H = E^{-1} = \begin{pmatrix} H' &  0 \\ H' & \alpha H' \end{pmatrix} \:,
\]
where again $H'$ is the analog for ${\bf s'}=(s,\ldots,s)$ (of length $r-1$) instead of~${\bf s}$.
The claim now follows by induction.
\\
(c) Combine (a) with Proposition~\ref{Prop:PDBasis}.
\end{proof}
\begin{cor} \label{Cor:EuT_HT_SL2x}
Let $x\in\{\id,s\}$.
\\
(a)
Let $\delta$, $\gamma\in\Gamma_x$. Then
\[
      \Eu_T(\delta,\ovl{C_x^\gamma})
    = \begin{cases} (-1)^{\#D(\delta)}(-\alpha)^{\# D(\gamma)}\,, &\text{if } D(\delta)\subseteq D(\gamma)  \\
                    \infty\,, &\text{else.}  \end{cases}
\]
(b)
The basis $\{\mu_{\gamma,x}\mid\gamma\in\Gamma_x\}$ of $\rH_T^\blt(\BS_x)$ is determined by
\[
    \Res{\mu_{\gamma,x}}{\delta}
  = \begin{cases} \alpha^{\#D(\gamma)}, &\text{if } D(\delta)\supseteq D(\gamma) , \\
                            0,          &\text{else.}                              \end{cases}
\]
(c)
The basis $\{\mu_{\gamma,x}^*\mid\gamma\in\Gamma_x\}$ of $\rH_T^\blt(\BS_x)$ is determined by
\[
    \Res{\mu_{\gamma,x}^*}{\delta}
  = \begin{cases} (-\alpha)^{r-1-\#D(\gamma)}, &\text{if } D(\delta)\subseteq D(\gamma) , \\
                                   0,          &\text{else.}                              \end{cases}
\]
\end{cor}
\begin{proof}
(a)
Use Equation~\eqref{Eq:EuTProd} and the fact that for $\pi(\delta)=\pi(\gamma)$ we have 
$J(\delta)\subseteq J(\gamma) \iff D(\delta)\subseteq D(\gamma)$.
\\
To deduce (b) and (c) we notice that $\iota(\BS_x)=(G/B)^{r-1}\times\{x\}$ and hence $\BS_x$ is isomorphic
to a Bott--Samelson variety $\BS'$ for  the sequence ${\bf s'}=(s,\ldots,s)$ (of length $r-1$) instead
of~${\bf s}$.
Under this isomorphism $\gamma\in\Gamma_x$ corresponds to $\gamma'\in\Gamma'$. Furthermore we have
$\#D(\gamma)=\#J(\gamma')$ and $D(\delta)\subseteq D(\gamma) \iff J(\delta')\subseteq J(\gamma')$.
Now the claims follow from the proposition.
\end{proof}
\begin{prop} \label{Prop:SL2Res}
(a)
For any $\mu=\sum\limits_{\gamma\in\Gamma} a_\gamma \, \mu_\gamma \in\rH_T^\blt(\BS)$ we have
\[
    a_\gamma = \frac{\sum\limits_{\ssstyle\delta\colon J(\delta)\subseteq J(\gamma)}
                       (-1)^{\#J(\delta)}\cdot\Res{\mu}{\delta}}{(-\alpha)^{\#J(\gamma)}}   \quad(\gamma\in\Gamma).
\]
In particular, $\ff \in \rH_T^\blt(\BS^T)$ is in the image of
the restriction $\rH_T^\blt(\BS)\tto \rH_T^\blt(\BS^T)$ if and only if
\begin{equation} \label{Eq:SL2Res}
     \sum_{\delta\in\Gamma\colon J(\delta)\subseteq J(\gamma)} (-1)^{\#J(\delta)} \, \ff(\delta)
                           \equiv 0  \pmod{\alpha^{\#J(\gamma)}}  
\end{equation}
for all $\gamma\in\Gamma$.
\\
(b)
For any $\mu=\sum\limits_{\gamma\in\Gamma} a_\gamma \, \mu_\gamma^* \in\rH_T^\blt(\BS)$ we have
\[
    a_\gamma = \frac{\sum\limits_{\ssstyle\delta\colon J(\delta)\supseteq J(\gamma)}
                       (-1)^{r-\#J(\delta)}\cdot\Res{\mu}{\delta}}{\alpha^{r-\#J(\gamma)}}   \quad(\gamma\in\Gamma).
\]
In particular, $\ff \in \rH_T^\blt(\BS^T)$ is in the image of
the restriction $\rH_T^\blt(\BS)\tto \rH_T^\blt(\BS^T)$ if and only if
\begin{equation} \label{Eq:SL2ResPD}
     \sum_{\delta\in\Gamma\colon J(\delta)\supseteq J(\gamma)} (-1)^{r-\#J(\delta)} \, \ff(\delta)
                           \equiv 0  \pmod{\alpha^{r-\#J(\gamma)}}  
\end{equation}
for all $\gamma\in\Gamma$.
\end{prop}
\begin{proof}
(a) follows from the fact that $E$ is the inverse of $H$ together with the explicit formula
in Proposition~\ref{Prop:EuT_HT_SL2}(a)
\\
(b)
Using the remark after Proposition~\ref{Prop:EuT_HT_SL2} we see that the inverse of $H^*$ is obtained from $E$
by reversing the order of the rows and columns (we assume they are ordered according to the lexicographic order)
and substituting $(-\alpha)$ for $\alpha$. Now we may argue as above. 
\end{proof}
Similarly we get:
\begin{prop} \label{Prop:SL2Resx}
(a)
For any $\mu=\sum\limits_{\gamma\in\Gamma_x} b_\gamma \, \mu_{\gamma,x} \in\rH_T^\blt(\BS_x)$ we have
\[
    b_\gamma = \frac{\sum\limits_{\ssstyle\delta\in\Gamma_x\colon D(\delta)\subseteq D(\gamma)}
                       (-1)^{\#D(\delta)}\cdot\Res{\mu}{\delta}}{(-\alpha)^{\#D(\gamma)}}   \quad(\gamma\in\Gamma_x).
\]
In particular, $\ff_x \in \rH_T^\blt(\Gamma_x)$ is in the image of
the restriction $\rH_T^\blt(\BS_x)\tto \rH_T^\blt(\Gamma_x)$ if and only if
\begin{equation} \label{Eq:SL2Resx}
     \sum_{\delta\in\Gamma_x\colon D(\delta)\subseteq D(\gamma)} (-1)^{\#D(\delta)} \, \ff_x(\delta)
                           \equiv 0  \pmod{\alpha^{\#D(\gamma)}}
\end{equation}
for all $\gamma\in\Gamma_x$. 
\\
(b)
For any $\mu=\sum\limits_{\gamma\in\Gamma_x} b_\gamma \, \mu_{\gamma,x}^* \in\rH_T^\blt(\BS_x)$ we have
\[
    b_\gamma = \frac{\sum\limits_{\ssstyle\delta\in\Gamma_x\colon D(\delta)\supseteq D(\gamma)}
                       (-1)^{r-1-\#D(\delta)}\cdot\Res{\mu}{\delta}}{\alpha^{r-1-\#D(\gamma)}}   \quad(\gamma\in\Gamma_x).
\]
In particular, $\ff_x \in \rH_T^\blt(\Gamma_x)$ is in the image of
the restriction $\rH_T^\blt(\BS_x)\tto \rH_T^\blt(\Gamma_x)$ if and only if
\begin{equation} \label{Eq:SL2ResxPD}
     \sum_{\delta\in\Gamma_x\colon D(\delta)\supseteq D(\gamma)} (-1)^{\#D(\delta)} \, \ff_x(\delta)
                           \equiv 0  \pmod{\alpha^{r-1-\#D(\gamma)}}
\end{equation}
for all $\gamma\in\Gamma_x$. 
\end{prop}
Note that $\Gamma_s=\{\gamma\in\Gamma\mid r\in J(\gamma)\}$ and
$\Gamma_\id=\{\gamma\in\Gamma\mid r\not\in J(\gamma)\}$.
Therefore replacing $\gamma_r$ by $s_{\alpha_r}\gamma_r$ defines an involution
(``folding the ends'') on $\Gamma$ that exchanges $\Gamma_s$ and $\Gamma_\id$.
Denote it by $\gamma\mapsto\ovl{\gamma}$.

Consider the restriction
\[
   \rH_T^\blt(\BS)\injto\rH_T^\blt(\BS_s\sqcup\BS_\id)=\rH_T^\blt(\BS_\id)\oplus\rH_T^\blt(\BS_s) \:.
\]
It is injective because the composition
$\rH_T^\blt(\BS)\to\rH_T^\blt(\BS_\id)\oplus\rH_T^\blt(\BS_s)\to\rH_T^\blt(\BS^T)$ is.
\begin{prop}
The pair $(\nu,\sigma)=\bigl(\sum_{\gamma\in\Gamma_\id} a_\gamma\,\mu_\id^\gamma,
\sum_{\gamma\in\Gamma_s} b_\gamma\,\mu_s^\gamma\bigr)$
is in the image of the restriction $\rH_T^\blt(\BS) \injto \rH_T^\blt(\BS_\id)\oplus\rH_T^\blt(\BS_s)$
if and only if $a_{\ovl{\gamma}}\equiv b_\gamma \pmod{\alpha}$ for all $\gamma\in\Gamma_s$.
\end{prop}
\begin{proof}
If $\gamma\in\Gamma_\id$ then \eqref{Eq:SL2Res} is equivalent to \eqref{Eq:SL2Resx}.
Therefore $(\nu,\sigma)$ is in the image of the restriction iff the conditions~\eqref{Eq:SL2Res}
are satisfied for all $\gamma\in\Gamma_s$, i.~e.\ if
\[
      \sum_{\ssstyle\delta\in\Gamma_\id\colon J(\delta)\subseteq J(\gamma)} (-1)^{\#J(\delta)}\,\ff(\delta)
    + \sum_{\ssstyle\delta\in\Gamma_\id\colon J(\delta)\subseteq J(\gamma)} (-1)^{\#J(\delta)}\,\ff(\delta)
                           \equiv 0  \pmod{\alpha^{\#J(\gamma)}} \:. 
\]
(Here $\ff$ is the restriction of $(\nu,\sigma)$ to $\rH_T^\blt(\BS^T)$.)
Since $\#J(\delta)=\#D(\delta)+1$, if $\delta\in\Gamma_s$, and $\#J(\delta)=\#D(\delta)$, else,
this condition is equivalent to
\[
            \frac{\sum\limits_{\ssstyle\delta\in\Gamma_\id\colon D(\delta)\subseteq D(\gamma)}
                       (-1)^{\#D(\delta)}\cdot \ff(\delta)}{(-\alpha)^{\#D(\gamma)}}
    \equiv  \frac{\sum\limits_{\ssstyle\delta\in\Gamma_s\colon D(\delta)\subseteq D(\gamma)}
                       (-1)^{\#D(\delta)}\cdot \ff(\delta)}{(-\alpha)^{\#D(\gamma)}}         \pmod{\alpha} \:. 
\]
for all $\gamma\in\Gamma_s$. The claim follows since $D(\gamma)=D(\ovl{\gamma})$.
\end{proof}
\begin{rem}
Comparing Proposition~\ref{Prop:EuT_HT_SL2}(b) to Corollary~\ref{Cor:EuT_HT_SL2x}(b)
one finds that $\mu_\gamma\mapsto(\mu_\id^\gamma,\mu_s^{\ovl{\gamma}})$,
if $\gamma\in\Gamma_\id$, and $\mu_\gamma\mapsto(0,\alpha\cdot\mu_s^\gamma)$,
if $\gamma\in\Gamma_s$. This gives another proof of the proposition.
Furthermore, it follows that the restrictions $\rH_T^\blt(\BS) \to \rH_T^\blt(\BS_x)$
are surjective.
\end{rem}
For $x\in\{\id,s\}$ denote by $B_x$ the image of the restriction map
$\rH_T^\blt(\BS_x)\to\rH_T^\blt(\Gamma_x)$.
(This is also the image of $\rH_T^\blt(\BS)\to\rH_T^\blt(\Gamma_x)$.)
Consider the $A$--linear map
$\rH_T^\blt(\BS_\id)\to \rH_T^\blt(\BS_s)$ sending the basis
element $\mu_\id^\gamma$ to $\mu_s^{\ovl{\gamma}}$ (for $\gamma\in\Gamma_\id$).
Again from Corollary~\ref{Cor:EuT_HT_SL2x}(b) we see that it induces an $A$--linear map
$\ovl{\rule{0pt}{1.4ex}\;\;}\colon B_\id \to B_s$, $\ff\mapsto\ovl{\ff}$ where 
$\ovl{\ff}(\delta)=\ff(\ovl{\delta})$.

Therefore defining
$\phi_\id\colon B_\id \to B_s/\alpha B_s$, ${\bf f}\mapsto {\bf \ovl{f}} + \alpha B_s$
and denoting
$\phi_s\colon B_s \to B_s/\alpha B_s$, ${\bf f}\mapsto {\bf f} + \alpha B_s$
the natural map we may reformulate the last proposition as follows:
\begin{cor} \label{Cor:SL2HTSheaf}
The image of the restriction map $\rH_T^\blt(\BS)\injto\rH_T^\blt(\BS^T)$
is contained in $B_\id\oplus B_s$ and $({\bf f}_\id,{\bf f}_s)\in B_\id\oplus B_s$
is in the image if and only if $\phi_\id({\bf f}_\id)=\phi_s({\bf f}_s)$.
\end{cor}
%
%
We conclude this section with a remark that will not be used in the rest of the paper.
There is an involution $\omega$ on the set of galleries which corresponds to reversing the
order of the factors of $(G/B)^r$:
Let $\omega(\delta)$ be defined by the requirement that $i\in J(\omega(\delta))$ if and
only if $r+1-i\in J(\delta)$.
\begin{lemma}
The matrices $E$ resp.\ $H$ are equivariant under the involution $\omega$ in the sense that
$\Eu_T(\omega(\delta),\ovl{C^{\omega(\gamma)}})=\Eu_T(\delta,\ovl{C^\gamma})$ resp.\
$\Res{\mu_{\omega(\delta)}}{\omega(\gamma)}=\Res{\mu_\delta}{\gamma}$.
\end{lemma}
\begin{proof}
Directly from Proposition~\ref{Prop:EuT_HT_SL2}(a) and (b).
\end{proof}
%
%
\section{The general Case}
\label{Sec:GeneralCase}
For any $T$--invariant subsets $Z_1\subseteq Z_2\subseteq\BS$ denote by
$\res_{Z_2,Z_1}\colon \rH_T^\blt(Z_2)\to\rH_T^\blt(Z_1)$ the restriction map. 
The following theorem is a special case of \cite[Theorem 6]{Br98}:
\begin{thm} \label{Thm:Brion}
The restriction map $\res_{\BS,\BS^T}\colon \rH_T^\blt(\BS)\to\rH_T^\blt(\BS^T)$ is injective
and its image is the intersection of the images of the restriction maps
$\res_{\BS^{T_\alpha},\BS^T}\colon \rH_T^\blt(\BS^{T_\alpha})\to\rH_T^\blt(\BS^T)$
over all $\alpha\in R$.
\end{thm}
\begin{rem}
In fact, to apply the theorem as stated in \cite{Br98} we have to replace $T$ by a compact
subtorus $K$ (which does not change the equivariant cohomology) and use the fact
that $\BS$ is a compact Hamiltonian $K$--space (cf. the remarks in~\cite{Br98}).
\\
An analog statement holds for $\BS_x$ instead of $\BS$.
\end{rem}
As an immediate application we get the following description of $\rH_T^\blt(\BS)$:
\begin{thm}\label{Thm:HTBS1}
Let $\ff\in\rH_T^\blt(\BS^T)$. Then the following are equivalent:
\begin{enumerate}
\item
${\bf f}$ is in the image of the restriction $\rH_T^\blt(\BS)\injto\rH_T^\blt(\BS^T)$.
\item
For all $\alpha\in R^+$ and all $\gamma\in\Gamma$ we have
\[
         \sum_{\substack{\delta\in\Gamma\colon\delta\sim_\alpha\gamma \\
                          J_\alpha(\delta)\subseteq J_\alpha(\gamma)}}
               (-1)^{\#J_\alpha(\delta)} \, \ff(\delta)
  \equiv 0
     \pmod{\alpha^{\#J_\alpha(\gamma)}} \:.
\]
\item
For all $\alpha\in R^+$ and all $\gamma\in\Gamma$ we have 
\[
         \sum_{\substack{\delta\in\Gamma\colon\delta\sim_\alpha\gamma \\
                          J_\alpha(\delta)\supseteq J_\alpha(\gamma)}}
               (-1)^{\#J_\alpha(\delta)} \, \ff(\delta)
  \equiv 0
     \pmod{\alpha^{\#M_\alpha(\gamma)-\#J_\alpha(\gamma)}} \:.
\]
\end{enumerate}
\end{thm}
\noindent
For any $x\in\pi(\Gamma)$ define
\[
   B_x^\alpha = \im\bigl(\res_{\BS_x^{T_\alpha},\Gamma_x}\colon
                           \rH_T^\blt(\BS_x^{T_\alpha})\to\rH_T^\blt(\Gamma_x)\bigr) \:.
\]
\begin{rems}
(1) The results of the last section show that $B_x^\alpha$ is a free $A$--module.
\\
(2) Because of the surjectivity of $\res_{\BS^{T_\alpha},\BS_x^{T_\alpha}}$ we also have
\[
   B_x^\alpha = \im\bigl(\res_{\BS^{T_\alpha},\Gamma_x}\colon
                           \rH_T^\blt(\BS^{T_\alpha})\to\rH_T^\blt(\Gamma_x)\bigr) \:.
\]
\end{rems}
\noindent
Now let
\[
   F_x = \bigcap\limits_{\alpha\in R^+} B_x^\alpha  \:.
\]
We will see shortly that $F_x$ is (an isomorphic image of) the $T$--equivariant cohomology
of the {\em fibre\/} $\BS_x$ (under the map $\res_{\BS_x,\Gamma_x}$).
\begin{lemma} \label{L:alphaFx}
Let $x\in\pi(\Gamma)$ and $\alpha\in R^+$.
\\
(a) Let $\ff\in\rH_T^\blt(\BS_x)$. Then the following are equivalent:
\begin{enumerate}
\item
$\ff \in B_x^\alpha$.
\item
For all $\gamma\in\Gamma_x$ we have
\begin{equation} \label{Eq:finBx}
       \sum_{\substack{\delta\in\Gamma_x\colon \delta\sim_\alpha\gamma\\
                       D_\alpha(\delta)\subseteq D_\alpha(\gamma)}}
        (-1)^{\#D_\alpha(\delta)} \, \ff(\delta)  \equiv 0  \pmod{\alpha^{\#D_\alpha(\gamma)}} \:.
\end{equation}
\item
For all $\gamma\in\Gamma_x$ we have
\[
       \sum_{\substack{\delta\in\Gamma_x\colon \delta\sim_\alpha\gamma\\
                       D_\alpha(\delta)\supseteq D_\alpha(\gamma)}}
        (-1)^{\#D_\alpha(\delta)} \, \ff(\delta)  \equiv 0  \pmod{\alpha^{\#M_\alpha(\gamma)-1-\#D_\alpha(\gamma)}} \:.
\]
\end{enumerate}
(b)
$\alpha F_x = F_x \cap \alpha B_x^\alpha$.
\end{lemma}
\begin{proof}
(a)
follows from Proposition~\ref{Prop:SL2Resx}.
\\
(b)
Clearly, $\alpha F_x \subseteq F_x \cap \alpha B_x^\alpha$.
So let $\ff \in \alpha B_x^\alpha$.
Then $\frac{1}{\alpha}\ff \in B_x^\alpha\subseteq \rH_T^\blt(\Gamma_x)$.
Now for any $\beta \in R^+$ different from $\alpha$ part (a) shows that
$\ff \in B_x^{\beta} \iff \frac{1}{\alpha}\ff \in B_x^{\beta} \iff \ff \in \alpha B_x^{\beta}$
whence the claim.
\end{proof}
\begin{lemma} \label{L:Fx}
For any $x\in\pi(\Gamma)$ we have
\[
   F_x = \im\bigl(\res_{\BS,\Gamma_x}\colon\rH_T^\blt(\BS)\to\rH_T^\blt(\Gamma_x)\bigr) \:.
\]
\end{lemma}
\begin{proof}
The inclusion $\supseteq$ is trivial because 
$\res_{\BS,\Gamma_x}=\res_{\BS^{T_\alpha},\Gamma_x}\circ\res_{\BS,\BS^{T_\alpha}}$
for all $\alpha\in R$.

We show the inclusion $\subseteq$ inductively:
Let us assume that any $\fg \in F_x$ such that $\fg(\delta)=0$ for all $\delta\leq\gamma$
is contained in $\im\bigl(\res_{\BS,\Gamma_x}\bigr)$.
Now take some $\ff \in F_x$ such that $\ff(\delta)=0$ for all $\delta<\gamma$,
hence in particular for all $\delta\sim_\alpha\gamma$, $\delta\not=\gamma$
such that $D_\alpha(\delta)\subseteq D_\alpha(\gamma)$.
By Lemma~\ref{L:alphaFx}(a) we know that $\ff(\gamma)$ is divisible by the product
of all $\alpha^{\#D_\alpha(\gamma)}$, i.~e.\ by $\Res{\mu_{\gamma,x}}{\gamma}$.

It follows that there is a $\fg$ as above such that the difference $\ff-\fg$
is a multiple of $\res_{\BS,\Gamma_x}(\mu_{\gamma,x})$.
Since $\fg \in\im\bigl(\res_{\BS,\Gamma_x}\bigr)$ by induction we conclude that $\ff \in\im\bigl(\res_{\BS,\Gamma_x}\bigr)$, too.
\end{proof}
\begin{cor}  \label{Cor:Fx}
(a) For any $x\in\pi(\Gamma)$ we have
\[
   F_x = \im\bigl(\res_{\BS_x,\Gamma_x}\colon\rH_T^\blt(\BS_x)\to\rH_T(\Gamma_x)\bigr) \:.
\]
In particular, $F_x$ is isomorphic to $\rH_T^\blt(\BS_x)$. Hence it is a free $A$--module.
\\
(b) The restriction map $\res_{\BS,\BS_x}\colon\rH_T^\blt(\BS)\to\rH_T^\blt(\BS_x)$
is surjective.
\end{cor}
\begin{proof}
(a) From the lemma we get
\[
   F_x =         \im\bigl(\res_{\BS,\Gamma_x}\bigr)
       =         \im\bigl(\res_{\BS_x,\Gamma_x}\circ\res_{\BS,\BS_x}\bigr)
       \subseteq \im\bigl(\res_{\BS_x,\Gamma_x}\bigr)
       \subseteq \bigcap\limits_{\alpha_\in R^+}
                        \im\bigl(\res_{\BS_x^{T_\alpha},\Gamma_x}\bigr)
       = F_x \:.
\]
So we have equality in all places. Now the first claim is immediate.
\\
(b) This follows from
$\im\bigl(\res_{\BS_x,\Gamma_x}\circ\res_{\BS,\BS_x}\bigr)
= \im\bigl(\res_{\BS_x,\Gamma_x}\bigr)$ together with the injectivity
of $\res_{\BS_x,\Gamma_x}$.
(See the remark after Theorem~\ref{Thm:Brion}.)
\end{proof}
\begin{rem}
There is an alternative proof:
In our terminology \cite[Satz 4.1]{So02} states that
$\res_{\BS,\Gamma_x}\bigl(\langle\mu_\gamma\mid\gamma\in\Gamma_x\rangle\bigr)
= d_x \cdot \im(\res_{\BS,\Gamma_x})$.
On the other hand Proposition~\ref{Prop:EuTRes} shows that
$\res_{\BS,\Gamma_x}\bigl(\langle\mu_\gamma\mid\gamma\in\Gamma_x\rangle\bigr)
= d_x \cdot \langle\mu_x^\gamma\mid\gamma\in\Gamma_x\rangle=d_x\cdot F_x$.
\end{rem}
As a consequence of Lemma~\ref{L:alphaFx}(b) the map
$j_x^\alpha\colon F_x/\alpha F_x\to B_x^\alpha/\alpha B_x^\alpha$,
${\bf f}+\alpha F_x\mapsto {\bf f}+\alpha B_x^\alpha$ is injective.
For $x\in\pi(\Gamma)$ and $\alpha\in R^+$ such that $s_\alpha x<x$ define
$\phi_{s_\alpha x}^\alpha\colon B_{s_\alpha x}^\alpha\to B_x^\alpha/\alpha B_x^\alpha$,
$\ff\mapsto {\bf \ovl{f}}+\alpha B_x^\alpha$
(where $\ovl{\rule{0pt}{1.4ex}\;\;}$ is the folding of ends corresponding
to the $\alpha$--wall).
\begin{lemma}
Let $x\in\pi(\Gamma)$ and $\alpha\in R^+$ such that $s_\alpha x<x$.
Then there is a (necessarily unique) $A$--linear map
$\rho_{s_\alpha x}^\alpha\colon F_{s_\alpha x}\to F_x/\alpha F_x$
such that $j_x^\alpha\circ\rho_{s_\alpha x}^\alpha=\phi_{s_\alpha x}^\alpha$ on $F_x$.
\end{lemma}
\begin{proof}
Let ${\bf f}_{s_\alpha x} \in F_{s_\alpha x}$.
By Lemma~\ref{L:Fx} there is a $\mu\in\rH_T^\blt(\BS)$ such that
${\bf f}_{s_\alpha x}=\res_{\BS,\Gamma_{s_\alpha x}}(\mu)$.
Let ${\bf f}_x=\res_{\BS,\Gamma_x}(\mu)\in F_x\subseteq B_x^\alpha$.
Then $({\bf f}_{s_\alpha x},{\bf f}_x)
\in\im\bigl(\res_{\BS,\Gamma_{s_\alpha x}\dotcup\Gamma_x}\bigr)$ implies that
$\ovl{{\bf f}_{s_\alpha x}}+\alpha B_x^\alpha = {\bf f}_x + \alpha B_x^\alpha$
(cf.\ Corollary~\ref{Cor:SL2HTSheaf}).
Therefore $\rho_{s_\alpha x}^\alpha({\bf f}_{s_\alpha x})={\bf f}_x$ does the
trick.
\end{proof}
For $x$ and $\alpha$ as in the lemma denote by
$\phi_x^\alpha\colon B_x^\alpha \surjto B_x^\alpha/\alpha B_x^\alpha$ resp.\
$\rho_x^\alpha\colon F_x\surjto F_x/\alpha F_x$ the natural maps.
\begin{thm} \label{Thm:HTBS2}
The image of the restriction map $\rH_T^\blt(\BS)\injto\rH_T^\blt(\BS^T)$
is contained in $\bigoplus_{x\in\pi(\Gamma)} F_x$
and $({\bf f}_x)_{x\in\pi(\Gamma)}\in \bigoplus_{x\in\pi(\Gamma)} F_x$ is in the image
if and only if for all $x\in\pi(\Gamma)$ and all $\alpha\in R^+$ such that $s_\alpha x<x$
we have $\rho_{s_\alpha x}^\alpha({\bf f}_{s_\alpha x})=\rho_x^\alpha({\bf f}_x)$.
\end{thm}
\begin{rem}
In Section~\ref{Sec:CombBases} we will give a recursive construction of bases for
the $F_x$.
\end{rem}
\begin{proof}
By Theorem~\ref{Thm:Brion} and Corollary~\ref{Cor:SL2HTSheaf} it suffices to show
that for ${\bf f}_{s_\alpha x}\in F_{s_\alpha x}$ and ${\bf f}_x\in F_x$ we have
\[
    \rho_{s_\alpha x}^\alpha({\bf f}_{s_\alpha x})=\rho_x^\alpha({\bf f}_x)
    \quad(\text{in } F_x/\alpha F_x)
    \quad\iff\quad
    \phi_{s_\alpha x}^\alpha({\bf f}_{s_\alpha x})=\phi_x^\alpha({\bf f}_x)
     \quad(\text{in } B_x^\alpha/\alpha B_x^\alpha) \:.
\]
But this is clear by the definitions.
\end{proof}
\begin{prop}  \label{Prop:Kernels}
Let $x\in\pi(\Gamma)$. Then:
\\
(a)
$\bigcap\limits_{\alpha\in R^+\colon s_\alpha x<x}\ker(\rho_x^\alpha)
= \res_{\BS,\Gamma_x}\bigl(\Span\{\mu_\gamma\mid\gamma\in\Gamma_x\}\bigr)
= d_x \cdot \res_{\BS_x,\Gamma_x}\bigl(\Span\{\mu_{\gamma,x}\mid\gamma\in\Gamma_x\}\bigr)$.
\\
(b) $\bigcap\limits_{\alpha\in R^+\colon s_\alpha x>x}\ker(\rho_x^\alpha)
= \res_{\BS,\Gamma_x}\bigl(\Span\{\mu_\gamma^*\mid\gamma\in\Gamma_x\}\bigr)$.
\end{prop}
\begin{proof}
Note that in either case we have $\ker(\rho_x^\alpha)=\ker(\phi_x^\alpha)$.
\\
(a)
See the remark after Corollary~\ref{Cor:Fx}.
\\
(b)
The inclusion $\supseteq$ is clear from Theorem~\ref{Thm:HTBS2}.

We show the inclusion $\subseteq$ by induction:
Let us assume that any $\fg \in \bigcap_{\alpha\colon s_\alpha x>x}\ker(\rho_x^\alpha)$
such that $\fg(\delta)=0$ for all $\delta\geq\gamma$ is contained in
$\res_{\BS,\Gamma_x}\bigl(\{\mu_\gamma^*\mid\gamma\in\Gamma^x\}\bigr)$.
Now take $\ff \in \bigcap\limits_{\alpha\colon s_\alpha x>x}\ker(\rho_x^\alpha)$
such that $\ff(\delta)=0$ for all $\delta>\gamma$, hence in particular for all
$\delta\sim_\alpha\gamma$, $\delta\not=\gamma$ such that $D_\alpha(\delta)\supseteq D_\alpha(\gamma)$.
Using Proposition~\ref{Prop:SL2Resx}(b) $\ff\in\ker(\phi_x^\alpha)$ implies that $\ff(\gamma)$ is
divisible by $\alpha^{\#M_\alpha(\gamma)-1-\#D_\alpha(\gamma)}$. So we get that $\ff(\gamma)$
is divisible by $\prod_{\alpha\colon s_\alpha x>x}\alpha^{\#M_\alpha(\gamma)\#J_\alpha(\gamma)}$
i.~e.\ by $\Res{\mu_\gamma^*}{\gamma}$ (see Proposition~\ref{Prop:PDBasis}).

It follows that there is a $\fg$ as above such that $\ff-\fg$ is a multiple of
$\res_{\BS,\Gamma_x}(\mu_\gamma^*)$ and hence
$\ff \in \res_{\BS,\Gamma_x}\bigl(\{\mu_\gamma^*\mid\gamma\in\Gamma^x\}\bigr)$ by induction.
\end{proof}
\begin{rem}
Let $F_x^*=\bigcap\limits_{\alpha\in R^+\colon s_\alpha x>x}\ker(\rho_x^\alpha)$.
Let $H_y$ be as in \eqref{Eq:DefHy} and let
\[
     H_y^*=\bigl(\Res{\mu^*_\gamma}{\delta}\bigr)_{\delta,\gamma\in\Gamma}
\]
be the matrix with columns the restrictions of the basis vectors $\mu^*_\gamma$
to the $T$--fixed points in $\BS_y$.
Then denoting $D_y=d_y\cdot D$ and $E_y=H_y^{-1}$ the matrix of inverse Euler classes
for $\BS_y$ we obtain as at the end of Section 3 the matrix
\[
    H_y^{-1} H_y^* = E_y\cdot D_y^{-1} \cdot {}^t\!E_y
\]
that expresses a basis of $F_y^*$ in terms of a basis of $F_y$.
Again this works for {\em any\/} basis of $F_y$ (in particular, for the combinatorial
bases to be constructed in Section~\ref{Sec:CombBases}).
\end{rem}
%
%
$\rH_T^\blt(\BS)$ is an $\rH_T^\blt(G/B)$--module via $\pi^*$.
By induction on the degree using the surjectivity of restriction to the stalks
one can show the following proposition.
\begin{prop}
As an $\rH_T^\blt(G/B)$--module $\rH_T^\blt(\BS)$ is generated by the set
$\{\mu_\gamma\mid\gamma\in\Gamma_\id\}$.
\end{prop}
%
%
\section{Sheaves on Bruhat graphs}
\label{Sec:Sheaves}
The aim of this section is to reformulate our results in the language of \cite{BM01}.
Let $Y$ be the {\em Bruhat graph\/} of $W$, i.~e.\ the oriented graph with vertices $W$ such that there is
an arrow $L$ from $x$ to $y$ whenever $y<x$ and $y=tx$ for some (not neccessarily simple) reflection $t\in W$.
In this case denote by $\alpha_L$ the root corresponding to the reflection $t$.
Let $Y^0$ resp.\ $Y^1$ be the vertices resp.\ arrows of $Y$.
For any vertex $x$ let $D_x$ be the set of arrows starting in $x$ and let $U_x$ be the set of arrows terminating
in $x$.
By abuse of notation we will also denote by $Y$ the disjoint union of $Y^0$ and $Y^1$ and we will view this
as a topological space in two different ways:

(I) For any $x\in Y^0$ resp.\ $L\in Y^1$ set $x^\circ=\{x\}\cup D_x\cup U_x$ resp.\ $L^\circ=\{L\}$.
These subsets of $Y$ form a basis for the open sets of the {\em ``standard topology''} on $Y$,
i.~e.\ $X\subseteq Y$ is open if and only if whenever it contains a vertex $x$ then it also contains all arrows incident to $x$.
Equivalently, $X\subseteq Y$ is closed, if and only if it is a subgraph.

(II) For any $X\subseteq Y$ set $X^\circ=\bigcup_{x\in X}x^\circ$.
Then a basis for the open sets of the {\em ``$B$--topology''} is given by the sets $\{y\in W\mid y\geq x\}^\circ$
for all  $x\in W$.
\begin{rem}
Let $\wtilde{Y}$ be the union of all $T$--orbits in $G/B$ of complex dimension at most one.
Then we get the topological space $Y$ as a quotient space of $\wtilde{Y}$ by the $T$--action
where we use the analytic topology on $G/B$ to get the standard topology on $Y$ and the topology
given by $B$--invariant open subsets on $G/B$ to get the $B$--topology on $Y$.
\end{rem}
Unless otherwise stated we will always consider $Y$ with its standard topology and we will use the terms
$B$--open, $B$--flabby, etc.\ to emphasize that we mean open, flabby, etc.\ with respect to the $B$--topology.
Since the $B$--topology is coarser than the standard topology, any sheaf with respect to the standard topology
is also a sheaf with respect to the $B$--topology. Therefore, the following lemma makes sense for both topologies.
\begin{lemma} \label{L:flabbiness}
Let $\cF$ be a sheaf on $Y$ with restriction maps $\rho_{U,V}$.
Then the following conditions are equivalent:
\begin{enumerate}
\item
$\cF$ is $B$--flabby.
\item
If $U\subset V$ are $B$--open subsets of $Y$ such that $V\setminus U$ contains only one vertex then
$\rho_{V,U}$ is surjective.
\item
If $U\subset V$ are $B$--open subsets of $Y$ such that $V\setminus U$ contains only one vertex then
$\ker(\rho_{V,U}) \subseteq \im(\rho_{Y,V})$.
\item
For any $y\in W$ and $V=\{y\in W\mid y\not< x\}^\circ$, $U=\{y\in W\mid y\not\leq x\}^\circ$ we have
$\ker(\rho_{V,U}) \subseteq \im(\rho_{Y,V})$.
\end{enumerate}
\end{lemma}
\noindent The proof is left to the reader.

A sheaf on the topological space $Y$ is completely determined by its sections over the basic sets $y^\circ$
for all $y\in Y$ together with the restriction maps between these sections.
As in \cite{BM01} we will henceforth consider only sheaves of $A$--modules such that for every $L\in Y^1$ the module
$\cF(L^\circ)$ is annihilated by $\alpha_L$.

Note that the module of sections over $y^\circ$ is at the same time the stalk $\cF_y$ at $y$.
So we will henceforth view a sheaf $\cF$ on $Y$ as a collection of $A$--modules $\cF_x$ ($x\in Y^0$),
$A/\alpha_L A$--modules $\cF_L$ ($L\in Y^1$) and $A$--linear maps $\rho_{x,L}$ for all pairs $(x,L)$
such that $x$ is incident to $L$.
An important example is the structure sheaf $\cA$ on $Y$ where all the stalks are free modules of rank
one and where the restriction maps are the canonical projections.

For any subset $X$ of $Y^0\cup Y^1$ (not necessarily open) a {\em section\/} of $\cF$ over $X$ is given by elements
$f_x\in\cF_x$ for $x\in X\cap Y^0$ and $f_L\in\cF_L$ for $L\in X\cap Y^1$ such that $\rho_{x,L}(f_x)=f_L$
whenever $\rho_{x,L}$ is defined.
Let $\cF(X)$ be the $A$--module of all sections over $X$.
Note that $\cA(X)$ is a ring and $\cF(X)$ is a module over $\cA(X)$.
Of course, $\cF(X)$ is the module of sections of the sheaf $\cF$ over $X$ in the usual sense if $X\subseteq Y$ is open.

Actually we have constructed a sheaf $\cBS$ on $Y$ with support $\pi(\Gamma)$ in the last section.
Set $\cBS_x=F_x$ for $x\in\pi(\Gamma)$ and $\cBS_x=(0)$ else.
Set $F_L=F_x/\alpha F_x$ if $L$ connects $x$ and $s_{\alpha_L}x$ where $s_{\alpha_L}x<x$.
Finally, set $\rho_{x,L}=\rho_x^{\alpha_L}$ if both $x$ and $s_{\alpha_L}x$ belong to
$\pi(\Gamma)$ and $\rho_{x,L}=0$ else.
Then Theorem~\ref{Thm:HTBS2} may be rephrased as follows.
\begin{thm} \label{Thm:BSSheaf}
The equivariant cohomology of the Bott--Samelson variety $\BS$ is isomorphic to the global
sections of the sheaf $\cBS$ defined above.
\end{thm}
Let $Y_{\geq x}$ be the full subgraph of $Y$ of all vertices $y\geq x$.
Braden and Macpherson call the sheaf $\cF$ pure, if
\begin{enumerate}
\item[(P1)]
All modules $\cF_x$ are free.
\item[(P2)]
$\cF_L=\cF_x/\alpha_L\cF_x$ if $x$ is the starting point of $L$.
\item[(P3)]
For all $x$ the restrictions $\cF(Y_{\geq x}\setminus\{x\})\to\cF(U_x)$ and $\cF_x\to\cF(U_x)$ have the same image.
\end{enumerate}
Let us have a closer look at the third condition: Let $\cF$ be a sheaf on $Y$ that satisfies (P3).
Then any section of $\cF$ over $Y_{\geq x}\setminus\{x\}$ can be extended to a section over
$Y_{\geq x}$, i.~e.\ $\rho_{Y_{\geq x},Y_{\geq x}\setminus\{x\}}$ is a surjection.
Hence, using Lemma~\ref{L:flabbiness} we get
\begin{enumerate}
\item[(P3.a)]
$\cF$ is $B$--flabby.
\end{enumerate}
Moreover, any section of $\cF$ over $x^\circ$ can be extended to a section over
$Y_{\geq x}$, i.~e.\ $\rho_{Y_{\geq x},x^\circ}$ is a surjection, too.
Together with the $B$--flabbiness we get
\begin{enumerate}
\item[(P3.b)]
The restriction map $\cF(Y)\surjto\cF_x$ is surjective for any $x\in W$.
\end{enumerate}
Conversely, let $\cF$ be any sheaf on $Y$ that satisfies (P3.a) and (P3.b).
Then for all $x$ the restrictions $\rho_{Y,Y_{\geq x}\setminus\{x\}}$ and
$\rho_{Y,x^\circ}$ are surjective. Therefore, both $\rho_{Y_{\geq x}\setminus\{x\},U_x}$ and
$\rho_{x^\circ,U_x}$ have the same image as $\rho_{Y,U_x}$, i.~e.\ (P3) holds.
\begin{thm} \label{Thm:BSPurity}
The sheaf $\cBS$ is pure in the sense of Braden and Macpherson.
\end{thm}
\begin{proof}
We use the characterization given above.
Condition (P2) is clear by construction, (P1) and (P3.a) have been verified in Corollary~\ref{Cor:Fx}
and (P3.b) follows from Proposition~\ref{Prop:Kernels}(b).
\end{proof}
\begin{rem}
The fact that the $T$--equivariant cohomology of a Bott--Samelson variety equals the
global sections of a combincatorially pure sheaf can be deduced from the fact (whose
proof uses the decomposition theorem for perverse sheaves) that it is the direct sum
of $T$--equivariant intersection cohomology groups of some Schubert varieties
together with the results of \cite{BM01}.
As mentionned before, one motivation for the present article was to give an independent
proof that does not use the theory of perverse sheaves.
\end{rem}
As a direct consequence (see \cite{BM01}) we note
\begin{cor}
$\cBS$ is the direct sum of indecomposable pure sheaves.
\end{cor}
%
\section{Bases for the Cohomology of the Fibres}
\label{Sec:CombBases}
In this section we show how to construct recursively bases of $F_x$.
As before, let $\BS$ be the Bott--Samelson variety of the sequence
${\bf s}=(s_1,\ldots,s_r)$.
Denote by $\BS'$ the Bott--Samelson variety of the sequence ${\bf s'}=(s_1,\ldots,s_{r-1})$.
Similarly, let $\Gamma'$ be the set of all combinatorial galleries in $\BS'$
and let $\Gamma\to\Gamma'\colon\gamma\mapsto\gamma'=[\gamma_1,\ldots,\gamma_{r-1}]$
be the map that ``cuts off'' the last bend or crossing.
Let $\alpha=\beta_r$ be the last wall of some $\gamma\in\Gamma_x$.
Then the projection $\BS\surjto\BS'$ gives rise to a bijection between $\Gamma_x$
and $\Gamma'_x\dotcup\Gamma'_{s_\alpha x}$ and hence to an isomorphism between
$\rH_T(\Gamma_x)$ and $\rH_T({\Gamma'}_x)\oplus\rH_T({\Gamma'}_{s_\alpha x})$.
Let $({\ff'}_x,{\ff'}_{s_\alpha x})$ be the image of $\ff_x$ under this isomorphism.
\begin{lemma}
(a) Let $\alpha=\beta_r$. Then the above isomorphism identifies $B_x^\alpha$ with
the submodule of ${B'}_x^\alpha \oplus {B'}_{s_\alpha x}^\alpha$ consisting of all
$(\ff'_x,\ff'_{s_\alpha x})$ such that
${\phi'}_x^\alpha({\ff'}_x)={\phi'}_{s_\alpha x}^\alpha({\ff'}_{s_\alpha x})$.
\\
(b) Let $\beta\not=\alpha=\beta_r$. Then the above isomorphism identifies $B_x^\beta$
with ${B'}_x^\beta \oplus {B'}_{s_\alpha x}^\beta$.
\end{lemma}
\begin{proof}
For the proof we will assume that  $s_\alpha x<x$.
(Except for notational changes the proof for $s_\alpha x>x$ is exactly the same.)

(a)
If $r\not\in D(\gamma)$, i.~e.\ if the gallery $\gamma'$ ends in $s_\alpha x$
then the map $\delta\mapsto\delta'$ sends the set $\bigl\{\delta\in\Gamma_x \mid
\delta\sim_\alpha\gamma, D_\alpha(\delta)\subseteq D_\alpha(\gamma)\bigr\}$
bijectively to $\bigl\{\delta'\in\Gamma'_{s_\alpha x}\mid
\delta'\sim_\alpha\gamma',D_\alpha(\delta')\subseteq D_\alpha(\gamma')\bigr\}$.
Hence Lemma~\ref{L:alphaFx} shows
that $\ff_x\in B_x^\alpha$ implies ${\ff'}_{s_\alpha x}\in {B'}_{s_\alpha x}$
and, conversely, that ${\ff'}_{s_\alpha x}\in {B'}_{s_\alpha x}$ implies that
the conditions \eqref{Eq:finBx} in Lemma~\ref{L:alphaFx} are fulfilled for $\ff_x$ and
all $\gamma\in\Gamma_x$ such that $r\not\in D(\gamma)$.

If $r\in D(\gamma)$, i.~e.\ if $\gamma'$ ends in $x$
then $\delta\mapsto\delta'$ sends $\bigl\{\delta\in\Gamma_x \mid
\delta\sim_\alpha\gamma, D_\alpha(\delta)\subseteq D_\alpha(\gamma)\bigr\}$
bijectively to $\bigl\{\delta'\in\Gamma'_x\mid \delta'\sim_\alpha\gamma',
D_\alpha(\delta')\subseteq D_\alpha(\gamma')\bigr\} \dotcup
\bigl\{\delta'\in\Gamma'_{s_\alpha x}\mid \delta'\sim_\alpha\gamma',
D_\alpha(\delta')\subseteq D_\alpha(\gamma')\bigr\}$.

Now \eqref{Eq:finBx} can be written as
\begin{equation} \label{Eq:recFx1}
      \sum_{\ssstyle\substack{\delta'\in\Gamma'_x\colon \delta'\sim_\alpha\gamma' \\
                      D_\alpha(\delta')\subseteq D_\alpha(\gamma')}} (-1)^{\#D_\alpha(\delta')}\,\ff'_x(\delta')
    + \sum_{\ssstyle\substack{\delta'\in\Gamma_{s_\alpha x}\colon \delta'\sim_\alpha\gamma', \\
                      D_\alpha(\delta')\subseteq D_\alpha(\gamma')}} (-1)^{\#D_\alpha(\delta')}\,\ff'_{s_\alpha x}(\delta')
      \equiv 0  \pmod{\alpha^{\#D_\alpha(\gamma)}} \:. 
\end{equation}
Here the second summand does not change if we replace $\gamma'$ by $\ovl{\gamma'}\in\Gamma_{s_\alpha x}$
(folding the ''end'' along the $\alpha$--wall).
So we already know that if $\ff_x\in B_x$ then the second summand is divisible
by $\alpha^{\#D_\alpha(\gamma')}$, hence so is the first.
Using Lemma~\ref{L:alphaFx} again, we conclude that ${\ff'}_x\in {B'}_x$ and, furthermore,
\begin{equation}
            \frac{\sum\limits_{\ssstyle\substack{\delta'\in\Gamma'_x\colon \delta'\sim_\alpha\gamma' \\
                              D_\alpha(\delta')\subseteq D_\alpha(\gamma')}}
                       (-1)^{\#D_\alpha(\delta')}\cdot \ff'_x(\delta')}{(-\alpha)^{\#D_\alpha(\gamma')}}
    \equiv  \frac{\sum\limits_{\ssstyle\substack{\delta'\in\Gamma_{s_\alpha x}\colon \delta'\sim_\alpha\gamma', \\
                              D_\alpha(\delta')\subseteq D_\alpha(\ovl{\gamma'})}}
                       (-1)^{\#D_\alpha(\delta)}\cdot \ff'_{s_\alpha x}(\delta)}{(-\alpha)^{\#D_\alpha(\gamma')}}
    \pmod{\alpha} \:. 
\end{equation}
By Proposition~\ref{Prop:SL2Resx}(a) this implies 
${\phi'}_x^\alpha({\ff'}_x)={\phi'}_{s_\alpha x}^\alpha({\ff'}_{s_\alpha x})$.
Conversely, if ${\ff'}_x\in {B'}_x$ and ${\ff'}_{s_\alpha x}\in {B'}_{s_\alpha x}$ then both sides of the last
equation are in $A$. So ${\phi'}_x^\alpha({\ff'}_x)={\phi'}_{s_\alpha x}^\alpha({\ff'}_{s_\alpha x})$ implies
\eqref{Eq:recFx1},
i.~e. the conditions \eqref{Eq:finBx} in Lemma~\ref{L:alphaFx} are fulfilled for $\ff_x$ and
all $\gamma\in\Gamma_x$ such that $r\not\in D(\gamma)$.

The proof of (b) is similar (but easier).
\end{proof}
As an immediate consequence we get
\begin{cor}
The isomorphism
$\rH_T^\blt(\Gamma_x) \isoto \rH_T^\blt({\Gamma'}_x)\oplus\rH_T^\blt({\Gamma'}_{s_\alpha x})$
identifies $F_x$ with the submodule of ${F'}_x \oplus {F'}_{s_\alpha x}$
consisting of all $({\ff'}_x,{\ff'}_{s_\alpha x})$ such that
${\rho'}_x^\alpha({\ff'}_x)={\rho'}_{s_\alpha x}^\alpha({\ff'}_{s_\alpha x})$.
\end{cor}
This corollary permits to calculate recursively bases of the modules $F_x$
that are indexed by the galleries in $\Gamma_x$ in such a way that the basis element
indexed by $\gamma$ is supported on the set $\{\delta\in\Gamma_x\mid\delta>\gamma\}$:

Let us suppose again, that $s_\alpha x<x$.
(Up to notation the other case works the same.)
Suppose we have already found bases
$\{b'_{\gamma,s_\alpha x}\mid\gamma\in{\Gamma'}_{s_\alpha x}\}$ of ${F'}_{s_\alpha x}$ and
$\{b'_{\gamma,x}\mid\gamma\in{\Gamma'}_x\}$ of ${F'}_x$.
Then the pairs $(c'_\gamma,b'_{\gamma,s_\alpha x})\in {F'}_x \oplus {F'}_{s_\alpha x}$
where $c'_\gamma$ is any lift of ${\rho'}_{s_\alpha x}^\alpha(b'_{\gamma,s_\alpha x})$
to ${F'}_x$
together with the pairs $(\alpha\cdot b'_{\gamma,x},0)\in {F'}_x \oplus {F'}_{s_\alpha x}$
form a basis of $F_x\injto {F'}_x \oplus {F'}_{s_\alpha x}$.
\begin{rem}
The support property of our 'combinatorial' basis of the cohomology of the fibre
$\BS_x$ is shared by the 'geometric' basis
$\{\mu_{\gamma,x}\mid\gamma\in\Gamma_x\}$ derived from the Bialynicki--Birula
decomposition of the fibre.
This means in particular that the base change matrix between the combinatorial and
the geometric basis will be lower triangular (provided the bases are ordered using $<$).
\end{rem}
%
%

\end{document}